\documentclass[a4paper,12pt,reqno]{amsart}

\usepackage{amsmath}
\usepackage{epsfig}
\usepackage{color} %% urs}
\makeatletter

\newtheorem{theorem}{Theorem}
\newtheorem{proposition}[theorem]{Proposition}
\newtheorem{corollary}[theorem]{Corollary}
\newtheorem{remark}[theorem]{Remark}

\newcommand{\diag}{\mathop{\mathrm{diag}}}
\newcommand{\Sym}{\mathop{\mathrm{Sym}}}
\def\PerfProof{{\it Proof.\ }}

\makeatother

\begin{document}

\date{}

 %%\newpage
 \thispagestyle{empty}
 \pagestyle{empty}

     \begin{color}{blue}\LARGE\boldmath
     \begin{center}
     ~\\ \vspace{15mm}
      {\bf Coxeter Transformations, the McKay correspondence, and the Slodowy correspondence}
      \end{center}
      \Large
      \end{color}
      \unboldmath \vspace{2mm}

      \begin{color}{blue}
       \begin{center}\LARGE\boldmath
       ~\\  \vspace{5mm}
         Rafael Stekolshchik
      \end{center}
      \end{color}

      \begin{color}{blue}
       \begin{center}\LARGE\boldmath
       ~\\  \vspace{45mm}
        BIRS, 08w5060 Workshop \\
        {\bf \lq\lq Spectral Methods in Representation Theory of Algebras and
         Applications to the Study of Rings of Singularities\rq\rq } \\
        ~\\  \vspace{2mm}
         September 8, 2008
      \end{center}
      \end{color}

 %%\maketitle

 \Large
 \newpage
 ~\\ \vspace{15mm}

     \begin{color}{red}\LARGE\boldmath
     \begin{center}
     ~\\ \vspace{15mm}
      {\bf Abstract}
      \end{center}
      \Large
      \end{color}
      \vspace{2mm}

%% \DL{horosho by skazat' kak mozhno skoree o KAKOM iz kuchi raznyh rydov
%% H-P idet rech'. Ibo vyklyuchnaya f-la nizhe --- yavnaya f-la dlya ryada, a kto PO
%% OPREDELENIYU koeff. pri $t^n$?}

   In \cite{Ebl02}, Ebeling established  a connection between
   certain Poincar\'{e} series, the Coxeter transformation \begin{color}{blue}\Large{\bf C}\end{color},
   and the corresponding affine
   Coxeter transformation  \begin{color}{blue}\Large${\bf C}_a$\end{color}
   (in the context of the McKay correspondence).
   We consider the  generalized Poincar\'{e} series
   \begin{color}{blue}\Large\boldmath$[\widetilde{P}_G(t)]_0$\end{color}
   for the case of multiply-laced diagrams
   (in the context of the McKay-Slodowy correspondence)
   and extend the Ebeling theorem for this case:
 \begin{color}{blue}\Large\boldmath
 \begin{equation*}
       [\widetilde{P}_G(t)]_0 =
       \frac{\mathcal{X}(t^2)}{\tilde{\mathcal{X}}(t^2)},
 \end{equation*}
 \end{color} \vspace{2mm}\\
where \begin{color}{blue}\Large\boldmath$\mathcal{X}$\end{color} is
the characteristic polynomial of the Coxeter transformation and
\begin{color}{blue}\Large\boldmath$\tilde{\mathcal{X}}$\end{color}
 is the characteristic polynomial of the
corresponding affine Coxeter transformation.

We obtain that Poincar\'{e} series coincide for pairs of diagrams
obtained by folding:
 \begin{color}{blue}\Large\boldmath
 \begin{equation*}
       \frac{\mathcal{X}(\Gamma)}{\mathcal{X}(\tilde\Gamma)} =
       \frac{\mathcal{X}(\Gamma^f)}{\mathcal{X}(\tilde\Gamma^f)},
 \end{equation*}
 \end{color} \vspace{2mm}\\
 where \begin{color}{blue}\Large\boldmath$\Gamma$\end{color} is any
 ($A$, $D$, $E$ type) Dynkin diagram,
 \begin{color}{blue}\Large\boldmath$\tilde\Gamma$\end{color} is the extended Dynkin diagram,
 and the diagrams \begin{color}{blue}\Large\boldmath$\Gamma^f$\end{color} and
 \begin{color}{blue}\Large\boldmath$\tilde\Gamma^f$\end{color} are obtained by folding from
 \begin{color}{blue}\Large\boldmath$\Gamma$\end{color} and \begin{color}{blue}\Large\boldmath$\tilde\Gamma$\end{color},
 respectively.

\newpage
~\\
\begin{color}{blue}\Large\boldmath
\tableofcontents
\newpage
\end{color}

\vspace{5mm}
%%%%%%%%%%%%%%%%%%%%%% PLAIN PAGESTYLE
 \pagestyle{plain}
 \Large
\newpage
~\\ \vspace{3mm}
    \begin{color}{red}\LARGE\boldmath
      \section{\sc\bf The Coxeter transformation \\ (A bit of history)}
    \end{color}
 \vspace{3mm}

Given a root system
\begin{color}{blue}\boldmath$\varDelta$\end{color}, a
\begin{color}{blue}{\bf Coxeter transformation}\end{color} (or
\begin{color}{blue}{\bf Coxeter element}\end{color})
\begin{color}{blue}\boldmath$C$\end{color} is defined as
the product of all the reflections in the simple roots. (We are
speaking here only about diagrams which are trees). Notations:

 \begin{color}{blue}\boldmath$h$\end{color} is
 the order of the Coxeter transformation (\begin{color}{blue}{\bf Coxeter
 number}\end{color}),

 \begin{color}{blue}\boldmath$|\varDelta|$\end{color} is the number of roots in the
 root system
\begin{color}{blue}\boldmath$\varDelta$\end{color},

\begin{color}{blue}\boldmath$l$\end{color} is the number of
eigenvalues of the Coxeter transformation, i.e., the number of
vertices in the Dynkin diagram.

We have:
\begin{color}{blue}\LARGE\boldmath
\begin{equation*}
 \label{coxeter_result_1}
               hl = |\varDelta|,
\end{equation*}
\end{color}

(Coxeter, \cite{Cox51}; Kostant \cite{Kos59}). %% \DL
Let
\begin{color}{blue} \boldmath$m_i$\end{color} be the exponents of the eigenvalues of
\begin{color}{blue}\boldmath$C$\end{color}, (all the eigenvalues in the case considered
 here are of the form $e^{2\pi im_j/h}$),
\begin{color}{blue}\boldmath$|W|$\end{color} be the order of the Weyl
group \begin{color}{blue}\boldmath$W$\end{color}.

Then
\begin{color}{blue}\LARGE\boldmath
\begin{equation*}
%%\label{coxeter_result_2}
  |W|  = (m_1 + 1)(m_2 + 1)...(m_l + 1),
\end{equation*}
\end{color}

(Coxeter, \cite{Cox34} ; proved by Chevalley \cite{Ch55} and other
authors). Let
\begin{color}{blue}\boldmath$\varDelta_{+}\subset \varDelta$\end{color} be the subset of simple positive
roots \begin{color}{blue}\boldmath$\alpha_i \in
\varDelta_{+}$\end{color},

\begin{color}{blue}\boldmath$\beta = n_1\alpha_1 + \dots +
n_l\alpha_l$\end{color}  be the highest root in the root system
\begin{color}{blue}\boldmath$\varDelta$\end{color}. Then

\begin{color}{blue}\LARGE\boldmath
\begin{equation*}
%%\label{coxeter_result_2}
  h  = n_1 + n_2 + ... + n_l + 1.
\end{equation*}
\end{color}

(Coxeter \cite{Cox49}; Steinberg \cite{Stb59}).

\newpage
~\\ \vspace{3mm}
    \begin{color}{red}\LARGE\boldmath
      \section{\sc\bf The Coxeter transformation \\ (bicolored partition)}
    \end{color}
 \vspace{3mm}

A partition $S = S_1 \coprod S_2$ of the vertices of the graph
$\varGamma$ is said to be \begin{color}{blue}{\bf
bicolored}\end{color} if all edges of $\varGamma$ lead from $S_1$ to
$S_2$. (A bicolored partition exists for trees). The diagram
$\varGamma$ admitting a bicolored partition is said to be
 \begin{color}{blue}{\bf bipartite}\end{color}.
\vspace{2mm}\\

An orientation $\varLambda$ is said to be \underline{bicolored}, if
there is the corresponding \underline{sink-admissible sequence}.

\begin{color}{blue}\LARGE\boldmath
\begin{equation*}
   \{v_1,v_2,~...,~v_m,v_{m+1},v_{m+2},~...~v_{m+k}\}
\end{equation*}
\end{color}

of vertices in this orientation $\varLambda$, such that the
subsequences %%\DL{v 2 stroki:}
\begin{color}{blue}\LARGE\boldmath
\begin{equation*}
  \begin{split}
   S_1 = & \{v_1, v_2,~...,~v_m\}, \\
   S_2 = & \{v_{m+1},v_{m+2},~...,~v_{m+k}\}
  \end{split}
\end{equation*}
\end{color}
form a bicolored partition, i.e., all arrows go from $S_1$ to $S_2$.
%%\DL
The product $w_i\in W(S_i)$ of all generators of $W(S_i)$ is an
involution  for $i = 1,2$, i.e.,

\begin{color}{blue}\LARGE\boldmath
\begin{equation}
 \label{C_decomp}
    w_1^2 = 1,  \qquad w_2^2 = 1, \qquad C = w_1 w_2.
\end{equation}
\end{color}

For the first time (as far as I know), the technique of bipartite
graphs was used by R.~Steinberg, \cite{Stb59}.

\newpage
~\\ \vspace{3mm}
    \begin{color}{red}\LARGE\boldmath
      \section{\sc\bf The Cartan matrix (Generalized) }
    \end{color}
 \vspace{3mm}

{\LARGE
 The \begin{color}{blue}{\bf generalized}\end{color} Cartan matrix:}

\begin{color}{blue}\LARGE\boldmath
\begin{equation*}
 \begin{split}
 \label{Cartan_matr}
 & (C1) \hspace{3mm}  {k}_{ii} = 2 \text{ for } i=1,..,n, \\
 & \\
 & (C2)  \hspace{3mm}  -{k}_{ij} \in \mathbb{Z}_+ = \{0,1,2,...\}
      \text{ for } i \neq j, \\
 & \\
 & (C3)  \hspace{3mm} {k}_{ij} = 0  \hspace{2mm} \text{implies} \hspace{2mm}
       {k}_{ji} = 0 \text{ for } i,j=1,...,n.
 \end{split}
\end{equation*}
\end{color}
\vspace{5mm}\\

 A generalized Cartan matrix $K$ is said to be
 \begin{color}{blue}{\bf symmetrizable}\end{color} if
there exists an invertible diagonal matrix $U$ with positive integer
coefficients and a symmetric matrix {\bf B} such that
 $K = U${\bf B}.

 (see Moody, \cite{Mo68}; Kac \cite{Kac80}).
 %%\DL{i Kac i mudi
 %% povisli v vozduhe bez svyazi s pred. i posl}

\begin{color}{blue}\boldmath
\begin{equation*}
 K =
   \begin{cases}
        2{\bf B} & \text{ for $K$ symmetric } \\
        & \\
        U{\bf B} & \text{ for $K$ symmetrizable }
   \end{cases}
\end{equation*}
\end{color}\vspace{2mm}\\
where $U$ is a diagonal matrix, {\bf B} is a symmetric matrix.

\newpage
~\\ \vspace{3mm}
    \begin{color}{red}\LARGE\boldmath
      \section{\sc\bf The Cartan matrix (diagrams) }
    \end{color}
 \vspace{3mm}

The {\it diagram} $(\varGamma, d)$ is a finite set $\varGamma_1$ (of
edges) rigged with numbers $d_{ij}$ for all pairs $i, j\in\partial
\varGamma_1\subset \varGamma_0$ (vertices) in such a way that

\begin{color}{blue}\LARGE\boldmath
\begin{equation*}
 \begin{split}
 & (D1) \hspace{3mm}  {d}_{ii} = 2 \text{ for } i=1,..,n, \\
 & \\
 & (D2)  \hspace{3mm}  {d}_{ij} \in \mathbb{Z}_+ = \{0,1,2,...\}
      \text{ for } i \neq j, \\
 & \\
 & (D3)  \hspace{3mm} {d}_{ij} = 0  \hspace{2mm} \text{implies} \hspace{2mm}
       {d}_{ji} = 0 \text{ for } i,j=1,...,n.
 \end{split}
\end{equation*}
\end{color}

It is depicted by symbols

\begin{color}{blue}\LARGE\boldmath
\begin{displaymath}
   i\stackrel{(d_{ij}, d_{ji})}{\line(1,0){40}}j
\end{displaymath}
\end{color}

If $d_{ij} = d_{ji} = 1$:

\begin{color}{blue}\LARGE\boldmath
\begin{displaymath}
   i \hspace{1mm} \line(1,0){40} \hspace{1mm} j
\end{displaymath}
\end{color}

There is a one-to-one correspondence between diagrams and
generalized Cartan matrices, and

\begin{color}{blue}\LARGE\boldmath
\begin{equation*}
    d_{ij} = |{k_{ij}}| \text{ for } i \neq j,
\end{equation*}
\end{color}
\vspace{2mm}\\
where $k_{ij}$ are elements of the Cartan matrix.

\newpage
~\\ \vspace{1mm}
    \begin{color}{red}\LARGE\boldmath
      \section{\sc\bf The Cartan matrix (simply-laced case) }
    \end{color}
 \vspace{2mm}

The integers $d_{ij}$ of the diagram are called
\begin{color}{blue}{\bf weights}\end{color},
and the corresponding edges are called
\begin{color}{blue}{\bf weighted edges}\end{color}.

The following edge is not weighted:
\begin{color}{blue}\LARGE\boldmath
\begin{equation*}
  \label{dij_dji}
 d_{ij} = d_{ji} = 1,
 \end{equation*}
 \end{color}

A diagram is called \begin{color}{blue}{\bf simply-laced}\end{color}
(resp. \begin{color}{blue}{\bf multiply-laced}\end{color}) if it
does not contain (resp. contains) weighted edges.
\vspace{2mm}\\

In the simply-laced case (= the symmetric Cartan matrix), we have:
\begin{color}{blue}\LARGE\boldmath
\begin{equation}
  \label{symmetric_B}
  \begin{split}
& K = 2{\bf B}, \text{\hspace{3mm}where\hspace{3mm}} {\bf B} = \left
(
\begin{array}{cc}
    I_m & D   \\
    D^t & I_k
\end{array}
\right ),  \\
 & \\
 & w_1 = \left (
\begin{array}{cc}
    -I_m & -2D   \\
    0 & I_k
\end{array}
\right ), \hspace{5mm} w_2 = \left (
\begin{array}{cc}
    I_m & 0   \\
    -2D^t & -I_k
\end{array}
\right ),
\end{split}
\end{equation}
\end{color} \\
where the elements $d_{ij}$ that constitute matrix $D$ are given by
the formula

\begin{color}{blue}\LARGE\boldmath
\begin{equation*}
  \label{matrix_D}
  d_{ij} = (a_i, b_j) =
  \left \{
   \begin{array}{cc}
     -\displaystyle\frac{1}{2} & \text{ if } |v(a_i) - v(b_j)| = 1 \hspace{1mm},
\vspace{2mm} \\
          0       & \text{ if } |v(a_i) - v(b_j)| > 1
          \hspace{1mm},
   \end{array}
  \right .
\end{equation*}
\end{color} \\
 where $v(a_i)$ and $v(b_j)$ are vertices lying in the
 different sets of the bicolored partition.

 \newpage
    \begin{color}{red}\LARGE\boldmath
      \section{\sc\bf The Cartan matrix (multiply-laced case) }
    \end{color}
 \vspace{3mm}

The multiply-laced case (= the symmetrizable and non-symmetric
Cartan matrix $K$):

\begin{color}{blue}\LARGE\boldmath
\begin{equation}
 \begin{split}
  \label{matrix_K}
& K = U{\bf B}, \hspace{3mm}\text{ where } \qquad
 K = \left (
\begin{array}{cc}
    2I_m & 2D   \\
    2F & 2I_k
\end{array}
\right ), \\
& w_1 = \left (
\begin{array}{cc}
    -I_m & -2D   \\
    0 & I_k
\end{array}
\right ), \hspace{4mm} w_2 = \left (
\begin{array}{cc}
    I_m & 0   \\
    -2F & -I_k
\end{array}
\right )
\end{split}
\end{equation}
\end{color} \\
with

\begin{color}{blue}\LARGE\boldmath
\begin{equation*}
     d_{ij} = \frac{(a_i, b_j)}{(a_i, a_i)},
     \hspace{10mm}
     f_{pq} = \frac{(b_p, a_q)}{(b_p, b_p)}, \vspace{3mm}
\end{equation*}
\end{color} \\
where the $a_i$ and $b_j$ are simple roots in the root systems
corresponding to $S_1$ and $S_2$, respectively. Here, $U = (u_{ij})$
is the diagonal matrix:
\begin{color}{blue}\Large\boldmath
\begin{equation*}
  u_{ii} = \frac{2}{(a_i, a_i)} = \frac{2}{\mathcal{B}(a_i)},
  \quad
   {\bf B} = \left ( \begin{array}{ccc}
       (a_i, a_i) & \dots & (a_i, b_j) \\
                  & \dots & \\
       (a_i, b_j) & \dots & (b_j, b_j) \\
  \end{array} \right ),
\end{equation*}
\end{color}

%%\begin{comment}
\begin{color}{blue}\LARGE\boldmath
\begin{equation*}
  K = U{\bf B} =  \left (
  \begin{array}{ccc}
       2 & \dots & \displaystyle\frac{2(a_i, b_j)}{(a_i, a_i)} \\
                  & \dots & \\
       \displaystyle\frac{2(a_i, b_j)}{(b_j, b_j)} & \dots & 2 \\
  \end{array} \right ).
\end{equation*}
\end{color}
%%\end{comment}

 \newpage
    \begin{color}{red}\LARGE\boldmath
      \section{\sc\bf The Cartan matrix (example: $\widetilde{F}_{41}$) }
    \end{color}
 \vspace{3mm}

The extended Dynkin diagrams $\widetilde{F}_{41}$ and
  $\widetilde{F}_{42}$

\begin{figure}[h]
\centering
\includegraphics{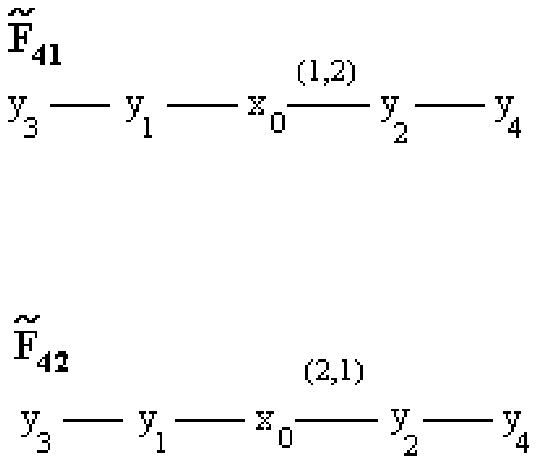}
\caption{\hspace{3mm} The diagrams $\widetilde{F}_{41}$ and
$\widetilde{F}_{42}$}
%%%%%% The label must come after caption
\label{F41_F42}
\end{figure}

a) Diagram $\widetilde{F}_{41}$. Here, the Cartan matrix is

\begin{color}{blue}\LARGE\boldmath
 \begin{equation*}
   \label{cartan_F41}
  K =
  \left (
   \begin{array}{ccccc}
     2 & -1 & -2 &    &       \\
    -1 & 2 &     & -1 &       \\
    -1 &   &  2 &     & -1    \\
       & -1  &    & 2 &      \\
       &   & -1   &   & 2    \\
   \end{array}
  \right )
   \begin{array}{c}
     x_0 \\
     y_1 \\
     y_2 \\
     y_3 \\
     y_4 \\
   \end{array}
 \end{equation*}
 \end{color}

 The matrix $U$ and the matrix {\bf B} of the
 \begin{color}{blue}{\bf Tits form}\end{color} are as follows:

\begin{color}{blue}\Large\boldmath
 \begin{equation*}
  \label{diag_F41}
  U = \diag \left (
   \begin{array}{c}
    1 \\
    1 \\
    1/2 \\
    1 \\
    1/2 \\
   \end{array}
  \right ),
  \qquad
  {\bf B} =
  \left (
   \begin{array}{ccccc}
     2 & -1 & -2 &    &       \\
    -1 & 2 &     & -1 &       \\
    -2 &   &  4 &     & -2    \\
       & -1  &    & 2 &       \\
       &   & -2   &   & 4     \\
   \end{array}
  \right ).
 \end{equation*}
 \end{color}

\newpage
    \begin{color}{red}\LARGE\boldmath
      \section{\sc\bf The Cartan matrix (example: $\widetilde{F}_{42}$) }
    \end{color}
 \vspace{3mm}

 b) Diagram $\widetilde{F}_{42}$. The Cartan matrix is
 \vspace{3mm}

\begin{color}{blue}\LARGE\boldmath
 \begin{equation*}
  \label{cartan_F42}
  K =
  \left (
   \begin{array}{ccccc}
     2 & -1 & -1 &    &       \\
    -1 & 2 &     & -1 &       \\
    -2 &   &  2 &     & -1    \\
       & -1  &    & 2 &      \\
       &   & -1   &   & 2    \\
   \end{array}
  \right )
   \begin{array}{c}
     x_0 \\
     y_1 \\
     y_2 \\
     y_3 \\
     y_4 \\
   \end{array},
 \end{equation*}
 \end{color} \\

 the matrix $U$ and the matrix {\bf B} of the Tits form
 are as follows:

\begin{color}{blue}\LARGE\boldmath
 \begin{equation*}
  \label{diag_F42}
  U = \diag \left (
   \begin{array}{c}
    1 \\
    1 \\
    2 \\
    1 \\
    2 \\
   \end{array}
  \right ),
  \qquad
  {\bf B} =
  \left (
   \begin{array}{ccccc}
     2 & -1 & -1 &    &       \\
    -1 & 2 &     & -1 &       \\
    -1 &   &  1 &     & -\displaystyle\frac{1}{2}    \\
       & -1  &    & 2 &       \\
       &   & -\displaystyle\frac{1}{2}   &   & 1     \\
   \end{array}
  \right ).
 \end{equation*}
 \end{color}

 \newpage
    \begin{color}{red}\LARGE\boldmath
      \section{\sc\bf The Cartan matrix and the Coxeter
transformation}
    \end{color}
 \vspace{3mm}

  From (\ref{C_decomp}), (\ref{matrix_K}) we have:

\begin{color}{blue}\Large\boldmath
\begin{equation}
\label{C_and_B_DF}
   {\bf C}{z} = \lambda{z}
   \hspace{3mm}
   \Longleftrightarrow
   \hspace{3mm}
   \left \{
       \begin{array}{cc}
         \displaystyle\frac{\lambda+1}{2\lambda}x & = -Dy   \\
         & \\
         \displaystyle\frac{\lambda+1}{2}y & = -Fx
       \end{array}
   \right . ,
   \text{ where }
   z = \left (
       \begin{array}{c}
          x \\
          y
       \end{array}
   \right ).
\end{equation}
\end{color}

\begin{color}{blue}\Large\boldmath
\begin{equation}
  \label{DDt_DtD}
\left \{
  \begin{array}{c}
     DD^{t}x = \displaystyle\frac{(\lambda+1)^2}{4\lambda}x \vspace{3mm} \\
     D^{t}Dy = \displaystyle\frac{(\lambda+1)^2}{4\lambda}y
  \end{array}
 \right.  \hspace{9mm}
 \left \{
  \begin{array}{c}
     DFx = \displaystyle\frac{(\lambda+1)^2}{4\lambda}x \vspace{3mm} \\
     FDy = \displaystyle\frac{(\lambda+1)^2}{4\lambda}y
  \end{array}
 \right.
\end{equation}
\end{color}
\vspace{5mm}\\

%%\DL
\begin{proposition}
 1) The kernel of the matrix {\bf B} considered as the matrix of an operator acting
 in the space spanned by roots coincides
with the kernel of the Cartan matrix $K$ and coincides with the
space of fixed points of the Coxeter transformation
\vspace{2mm}
\begin{color}{blue}\LARGE\boldmath
\begin{equation*}
   \ker{K} = \ker{\bf B} =
     \{ z \mid {\bf C}z = z \}.
\end{equation*}
\end{color}

2) The space of fixed points of the matrix {\bf B} coincides with
the space of anti-fixed points of the Coxeter transformation
\vspace{2mm}
\begin{color}{blue}\LARGE\boldmath
\begin{equation*}
   \{ z \mid {\bf B}z = z \} =
   \{ z \mid {\bf C}z = -z \}.
\end{equation*}
\end{color}
\end{proposition}

 \newpage
 ~\\ \vspace{3mm}
    \begin{color}{red}\LARGE\boldmath
      \section{\sc\bf The eigenvalues of the matrices $DF$ and $FD$}
    \end{color}
 \vspace{3mm}

1) The matrices $DF$ and $FD$
  have \begin{color}{blue}{\bf the same non-zero eigenvalues}\end{color}
   with equal multiplicities.

2) The eigenvalues  $\varphi_i$ of the matrices $DF$ and $FD$ are
non-negative:

\begin{color}{blue}\LARGE\boldmath
\begin{equation*}
      \varphi_i \geq 0.
\end{equation*}
\end{color}

3) The corresponding eigenvalues $\lambda^{\varphi_i}_{1,2}$ of
 the Coxeter transformations are
\begin{color}{blue}\LARGE\boldmath
  \begin{equation}
    \label{Coxeter_eigenvalues}
        \lambda^{\varphi_i}_{1,2} =
        2\varphi_i - 1 \pm 2\sqrt{\varphi_i(\varphi_i-1)}.
  \end{equation}
\end{color}

The eigenvalues $\lambda^{\varphi_i}_{1,2}$ either lie
  on the unit circle or are real positive numbers. It the latter case
  $\lambda^{\varphi_i}_1$ and $\lambda^{\varphi_i}_2$ are mutually inverse:

\begin{color}{blue}\LARGE\boldmath
  \begin{equation*}
     \lambda^{\varphi_i}_1\lambda^{\varphi_i}_2 = 1.
  \end{equation*}
\end{color}

 \newpage
 ~\\ \vspace{3mm}
    \begin{color}{red}\LARGE\boldmath
      \section{\sc\bf An example: a simple star  $\ast_{k+1}$}
    \end{color}
 \vspace{3mm}

In the simply-laced case, the following relation holds:

\begin{color}{blue}\LARGE\boldmath
\begin{equation*}
  \label{dd_formula}
  \begin{split}
     4(D&D^t)_{ij} =
     4\sum\limits_{p=1}^k{(a_i,b_p)(b_p, a_j)}  = \\
     &
       \begin{cases}
        s_i  & \text{if }  i = j, \\
        1    & \text{if }  |v_i - v_j| = 2, \\
        0    & \text{if }  |v_i - v_j| > 2,
      \end{cases}
  \end{split}
\end{equation*}
\end{color}
where $s_i$  is the number of edges with the  vertex $v_i$.
 %%\DL{as an endpoint}

\begin{figure}[h]
\centering
\includegraphics{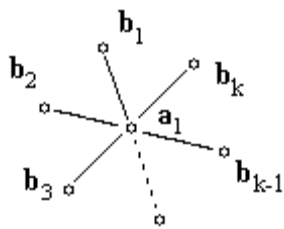}
\caption{\hspace{3mm} The star $*_{k+1}$ with $k$ rays}
%%%%%% The label must come after caption
\label{small_star}
\end{figure}

 \newpage
 ~\\ \vspace{3mm}
    \begin{color}{red}\LARGE\boldmath
      \section{\sc\bf An example: a simple star  $\ast_{k+1}$ (2) }
    \end{color}
 \vspace{3mm}

%%\DL{let?} n = k+1
 In the bicolored partition,
  one part of the graph consists of only one vertex $a_1$, i.e.,
  $m = 1$, the other one consists of $k$ vertices
  \{$b_1,\ldots,b_k$\}. Let $n = k+1$. The $1\times{1}$ matrix
  $DD^t$ is
\begin{color}{blue}\LARGE\boldmath
\begin{equation*}
        DD^t = k = n-1,
\end{equation*}
\end{color}
and the $k\times{k}$ matrix $D^tD$ is

\begin{color}{blue}\LARGE\boldmath
\begin{equation*}
       D^tD = \left (
        \begin{array}{ccccc}
          1 & 1 & 1 & \dots & 1 \\
          1 & 1 & 1 & \dots & 1 \\
          1 & 1 & 1 & \dots & 1 \\
            &   &   & \dots &   \\
          1 & 1 & 1 & \dots & 1
        \end{array}
        \right ).
\end{equation*}
\end{color}

The matrices $DD^t$ and $D^tD$ have only one non-zero eigenvalue
$\varphi_1 = n-1$. All the other eigenvalues of $D^tD$ are zeros and
the characteristic polynomial of the $D^tD$ is
\begin{color}{blue}\LARGE\boldmath
\begin{equation*}
         \varphi^{n-1}(\varphi - (n-1)).
\end{equation*}
\end{color}

 \newpage
    \begin{color}{red}\LARGE\boldmath
      \section{\sc\bf The Perron-Frobenius theorem}
    \end{color}
 \vspace{3mm}

%% \DL{kot takaya irr. matrix?!}
%% \DL{eto ne
%% t.e. a utv o tom, chto est' max. po abs. velichine i positive}

\begin{theorem}
   \label{pf_theorem}
 Let $A$ be an $n\times{n}$ non-negative irreducible matrix. Then the following holds:

 1) There exists a positive eigenvalue $\lambda$ such that
 \begin{color}{blue}\Large\boldmath
\begin{equation*}
   |\lambda_i| \leq \lambda, \text{ where } i = 1,2,\dots,n.
\end{equation*}
\end{color}

 2) There is a positive eigenvector $z$ corresponding to the eigenvalue $\lambda$:
\begin{color}{blue}\Large\boldmath
\begin{equation*}
   Az = \lambda{z}, \text{ where } z =(z_1,\dots,z_n)^t \text{ and }  z_i > 0
   \text { for } i = 1,\dots,n.
\end{equation*}
\end{color}
Such an eigenvalue $\lambda$ is called the
\begin{color}{blue}{\bf dominant eigenvalue}\end{color} of $A$.

 3) The eigenvalue $\lambda$ is a \underline{simple root} of the characteristic equation of $A$.
\end{theorem}

The eigenvalue $\lambda$ is calculated as follows:
\begin{equation*}
  \begin{split}
  & \displaystyle \lambda = \max_{z \geq 0} \min_i \frac{(Az)_i}{z_i}
     \hspace{7mm} (z_i \neq 0),       \\
  & \displaystyle \lambda = \min_{z \geq 0} \max_i \frac{(Az)_i}{z_i}
     \hspace{7mm} (z_i \neq 0).
  \end{split}
\end{equation*}

\newpage
    \begin{color}{red}\LARGE\boldmath
      \section{\sc\bf The Jordan normal forms of $DF$ and $FD$}
    \end{color}
 \vspace{3mm}

Here is an application of the Perron-Frobenius theorem.

\begin{figure}[h]
\centering
\includegraphics{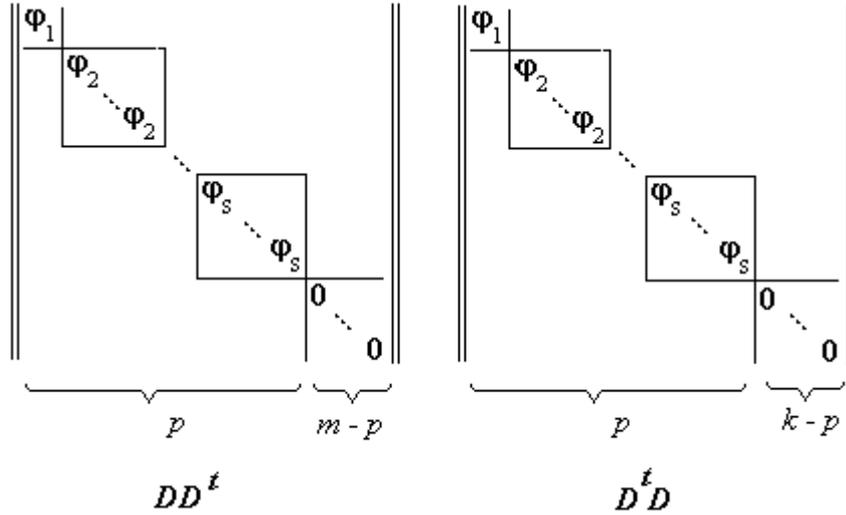}
\caption{\hspace{3mm} The Jordan normal forms of $DD^t$ and $D^tD$}
%%%%%% The label must come after caption
\label{normal_ddt}
\end{figure}

The matrices $DD^t$ (resp. $D^tD$) are symmetric and can be
diagonalized in the some orthonormal basis of the eigenvectors from
\begin{color}{blue}\Large\boldmath$\mathcal{E}_{\varGamma_a} = \mathbb{R}^h$\end{color} (resp.
\begin{color}{blue}\Large\boldmath$\mathcal{E}_{\varGamma_b} = \mathbb{R}^k$)\end{color}.
The Jordan normal forms
of these matrices are shown in Fig.~\ref{normal_ddt}.

%%\DL{vstav' slova:}
 In according to eq. (3.5), (3.14) from \cite{St08}, we have:
\begin{color}{blue}\Large\boldmath
\begin{equation*}
  U_1 A = D, \quad U_2 A^t = F, \quad DF = U_1 A U_2 A^t,
\end{equation*}
\end{color}\\
where $U_1$, $U_2$ are positive diagonal matrices, and the
eigenvalues of $DF$ and the symmetric matrix
\begin{color}{blue}\Large\boldmath
\begin{equation*}
   \sqrt{U_1}{A}U_2{A}^t\sqrt{U_1}
\end{equation*}
\end{color}\\
coincide.

The normal forms of $DF$ and $FD$ are the same, however, the normal
bases (i.e., bases which consist of eigenvectors) for $DF$ and $FD$
are not necessarily orthonormal: $\sqrt{U_1}$ does not preserve
orthogonality.

\newpage
    \begin{color}{red}\LARGE\boldmath
      \section{\sc\bf The eigenvectors of the Coxeter transformation}
    \end{color}
 \vspace{3mm}

\underline{Case $\varphi_i \neq 0,1$}:

\begin{color}{blue}\Large\boldmath
\begin{equation*}
 \label{case_non_01}
   z^{\varphi_i}_{r,\nu} =
   \left (
     \begin{array}{c}
        \mathbb{X}^{\varphi_i}_r  \vspace{3mm} \\
       -\displaystyle\frac{2}{\lambda^{\varphi_i}_{\nu}+1}
       D^t\mathbb{X}^{\varphi_i}_r
     \end{array}
   \right ), \hspace{4mm}
   1 \leq i \leq s, \hspace{2mm} 1 \leq r \leq t_i,
   \hspace{2mm} \nu = 1,2 \hspace{1mm}.
\end{equation*}
\end{color}\\

Here $\lambda^{\varphi_i}_{1,2}$ is obtained by eq.
(\ref{Coxeter_eigenvalues}). \\

\underline{Case $\varphi_i = 1$}:
\begin{color}{blue}\Large\boldmath
\begin{equation*}
 \label{case_1}
   z^1_{r} =
   \left (
     \begin{array}{c}
        \mathbb{X}^1_r  \vspace{3mm} \\
       -D^t\mathbb{X}^1_r
     \end{array}
   \right ), \hspace{4mm}
   \tilde{z}^1_{r} =
   \frac{1}{4}
   \left (
     \begin{array}{c}
        \mathbb{X}^1_r  \vspace{3mm} \\
       D^t\mathbb{X}^1_r
     \end{array}
   \right ), \hspace{4mm}
     1 \leq r \leq t_i.
\end{equation*}
\end{color}

\underline{Case $\varphi_i = 0$}:
\begin{color}{blue}\Large\boldmath
\begin{equation*}
 \label{case_0}
   z^0_{x_\eta} =
   \left (
     \begin{array}{c}
       \mathbb{X}^0_{\eta}  \vspace{3mm} \\
       0
     \end{array}
   \right ), \hspace{2mm}
   1 \leq \eta \leq m-p,
   \hspace{9mm}
   z^0_{y_\xi} =
   \left (
     \begin{array}{c}
       0      \vspace{3mm}     \\
       \mathbb{Y}^0_{\xi}
     \end{array}
   \right ), \hspace{2mm}
   1 \leq \xi \leq k-p.
\end{equation*}
\end{color}

These eigenvectors constitute the basis for the Jordan form of the
Coxeter transformation in the simply-laced case. (The multiply-laced
case is similarly considered, see \S3.2.2 and \S3.3.1 from
\cite{St08}.)
 %%\DL{:?}

\begin{color}{blue}\Large\boldmath
\begin{equation*}
  \label{PK_not_01}
   {\bf C}{z}^{\varphi_i}_{r,\nu} =
      \lambda^{\varphi_i}_{1,2}{z}^{\varphi_i}_{r,\nu}, \qquad
     \varphi_i \neq 0, 1.
\end{equation*}
\end{color}

\begin{color}{blue}\Large\boldmath
\begin{equation*}
  \label{PK_1}
   {\bf C}z^1_r = z^1_r, \qquad
   {\bf C}\tilde{z}^1_r = z^1_r + \tilde{z}^1_r,
   \qquad \varphi_i = 1, \quad \lambda = 1.
\end{equation*}
\end{color}

\begin{color}{blue}\Large\boldmath
\begin{equation*}
  \label{PK_0}
   {\bf C}z^0_{x_\eta} = -z^0_{x_\eta}, \qquad
   {\bf C}z^0_{y_\xi} = -z^0_{y_\xi},
   \qquad \varphi_i = 0, \quad \lambda = -1.
\end{equation*}
\end{color}

%%\newpage
    \begin{color}{red}\LARGE\boldmath
      \section{\sc\bf The Jordan form of the Coxeter transformation}
    \end{color}
 \vspace{3mm}

%%"if the Tits form"\DL{ne byla opredelena}.
%% R.S. Eto presentation, a polnyj tekst v knige
%% Opredelenie formy Titsa vstavit' v v.2
\begin{color}{blue}\Large\boldmath
\begin{theorem}
\label{th_jordan}
 1) The Jordan form of the Coxeter transformation is diagonal if and only if
  the Tits form is non-degenerate. \vspace{2mm}\\

 2) If $\mathcal{B}$ is non-negative definite ($\varGamma$ is an extended Dynkin diagram), then
  the Jordan form of the Coxeter transformation contains
  one $2\times2$ Jordan block. The remaining Jordan blocks are $1\times1$.
  All eigenvalues $\lambda_i$ lie on the unit circle. \vspace{2mm}\\

 3) If $\mathcal{B}$ is indefinite and degenerate, then
  the number of $2\times2$ Jordan blocks coincides with
  $\dim\ker{\bf B}$. The remaining Jordan blocks are $1\times1$.
  There is a simple maximal eigenvalue $\lambda^{\varphi_1}_1$
  and a simple minimal eigenvalue $\lambda^{\varphi_1}_2$, and
\begin{equation*}
   \lambda^{\varphi_1}_1 > 1, \hspace{5mm} \lambda^{\varphi_1}_2 < 1 .
\end{equation*}
\end{theorem}
\end{color}

  Subbotin-Stekolshchik, \cite{SuSt75}, \cite{SuSt78}.
 Similar results are obtained by A'Campo in \cite{A'C76}.

\begin{figure}[h]
\centering
\includegraphics{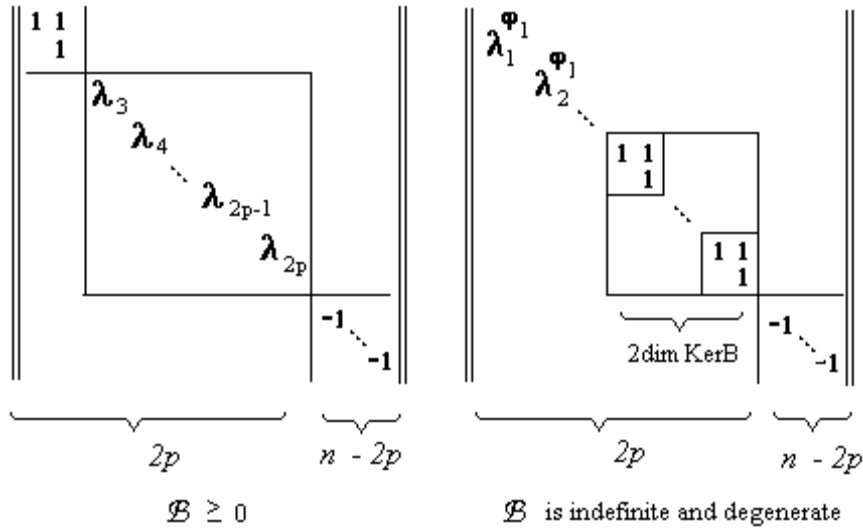}
\caption{\hspace{3mm} The Jordan normal form of the Coxeter
transformation}
%%%%%% The label must come after caption
\label{jordan_form}
\end{figure}

\newpage
    \begin{color}{red}\LARGE\boldmath
      \section{\sc\bf Example:
      an arbitrary large number of $2\times2$ Jordan blocks (Kolmykov)}
        \label{sect_blocks}
      \end{color}
    \vspace{3mm}

\begin{figure}[h]
\centering
\includegraphics{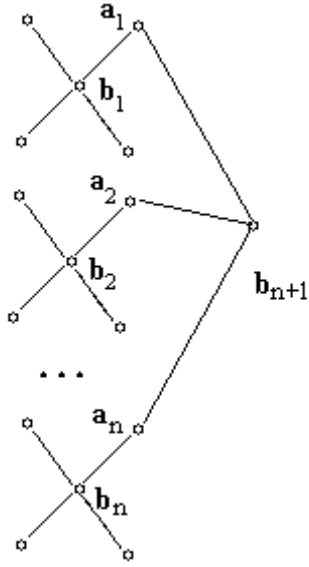}
\caption{\hspace{3mm} A graph $\varGamma$ such that $\dim\ker{\bf
B}$ is an arbitrary number}
%%%%%% The label must come after caption
 \label{dim_however_large}
\end{figure}

The example shows that there is a graph $\varGamma$ with indefinite
and degenerate quadratic form $\mathcal{B}$ such that $\dim\ker{\bf
B}$ is an arbitrarily large number (see Fig.
\ref{dim_however_large}) and the Coxeter transformation has an
arbitrary large number of $2\times2$ Jordan blocks.

%%\DL{slova:}
 We have:
\begin{color}{blue}\Large\boldmath
  \begin{equation*}
       4D^tD = \left (
        \begin{array}{ccccccc}
          n & 1 & 1 & 1 & \dots & 1 & 1 \\
          1 & 4 & 0 & 0 & \dots & 0 & 0 \\
          1 & 0 & 4 & 0 & \dots & 0 & 0 \\
          1 & 0 & 0 & 4 & \dots & 0 & 0 \\
            &   &   &  &  \dots &   &   \\
          1 & 0 & 0 & 0 & \dots & 4 & 0 \\
          1 & 0 & 0 & 0 & \dots & 0 & 4
        \end{array}
        \right ).
  \end{equation*}
 \end{color}
It is easy to show that
\begin{color}{blue}\Large\boldmath
  \begin{equation*}
    |4D^tD - \mu{I}| =
    (n - \mu)(4 - \mu)^n - n(4 - \mu)^{n-1} .
 \end{equation*}
 \end{color}
 Thus, \begin{color}{blue}\Large\boldmath$\varphi_i = \displaystyle\frac{\mu_i}{4} =
 1$\end{color}
  is of multiplicity
 \begin{color}{blue}\Large\boldmath\Large$n - 1$.\end{color}

\newpage
    \begin{color}{red}\LARGE\boldmath
      \section{\sc\bf Monotonicity of the dominant eigenvalue}
      \end{color}
    \vspace{3mm}

\begin{proposition}
Let us add an edge to a tree $\varGamma$ and let  %% \DL{kak ono pribavleno?!}
$\stackrel{\wedge}{\varGamma}$ be the new graph. Then:

  1) The dominant eigenvalue $\varphi_1$ may only grow:
\begin{color}{blue}\Large\boldmath
  \begin{equation}
   \label{only_grow}
   \varphi_1(\stackrel{\wedge}{\varGamma}) > \varphi_1(\varGamma)\hspace{1mm}.
  \end{equation}
\end{color}

  2) Let $\varGamma$ be an extended Dynkin diagram, i.e., $\mathcal{B}$
  is non-negative definite. Then the spectra of $DD^t(\stackrel{\wedge}{\varGamma})$ and
$D^tD(\stackrel{\wedge}{\varGamma})$ (resp.
$DF(\stackrel{\wedge}{\varGamma})$ and
 $FD(\stackrel{\wedge}{\varGamma})$) do not contain $1$, i.e.,
\begin{color}{blue}\Large\boldmath
  \begin{equation*}
   \varphi_i(\stackrel{\wedge}{\varGamma}) \neq 1
  \end{equation*}
\end{color}
for all $\varphi_i$ are eigenvalues
 of $DD^t(\stackrel{\wedge}{\varGamma})$.

 3) Let $\mathcal{B}$ be indefinite. Then
\begin{color}{blue}\Large\boldmath
  \begin{equation*}
      \varphi_1(\stackrel{\wedge}{\varGamma}) > 1\hspace{0.5mm}.
  \end{equation*}
\end{color}
\end{proposition}

  Subbotin-Stekolshchik, \cite{SuSt75}, \cite{SuSt78}.
  \vspace{2mm} \\

 During my talk Ringel noted  that
  (\ref{only_grow}) is a \underline{strict inequality}.
  The strict inequality (7)
  is, exactly, the
  result of \underline{Th. $1$} from \cite{SuSt78},
  and it is deduced from the following relation:

\begin{color}{blue}\Large\boldmath
 \begin{equation*}
  \label{3_char_pol}
  |DF(\stackrel{\wedge}{\varGamma}) - \mu{I}| =
    |DF - \mu{I}| +
    {\cos}^2\{a_i, b_s\}
     |DF(\stackrel{\vee}{\varGamma}) - \mu{I}|,
\end{equation*}
\end{color}\\
where
\begin{color}{blue}\Large\boldmath$\stackrel{\vee}{\varGamma}$\end{color}
is the diagram obtained from
\begin{color}{blue}\Large\boldmath$\varGamma$\end{color} by removing the vertex
$a_i$, and $b_s$ is the new vertex in the diagram
\begin{color}{blue}\Large\boldmath$\stackrel{\wedge}{\varGamma}$\end{color}.

\newpage
    \begin{color}{red}\LARGE\boldmath
      \section{\sc\bf Theorem on the spectral radius (Ringel)}
      \end{color}
    \vspace{3mm}

  The \begin{color}{blue}{\bf spectral radius}\end{color}
  $\rho(L)$ of a linear transforation $L$ of $\mathbb{R}^n$
  is the maximum of absolute values of the eigenvalues of $L$.
  The following theorem (due to C.~M.~Ringel \cite{Rin94})
  concerns the spectral radius of the Coxeter
  transformation in the case of the generalized Cartan matrix,
  including the case of \underline{diagrams with cycles}.
  \vspace{2mm}\\

\begin{color}{blue}\Large\boldmath
  \begin{theorem} Let $A$ be a generalized
  Cartan matrix which is connected and neither of finite
  nor of affine type. Let $C$ be a Coxeter transforation
  for $A$. Then $\rho(C) > 1$, and $\rho(C)$ is an
  eigenvalue of multiplicity  one, whereas any other
  eigenvalue $\lambda$ of $C$ satisfies
  $|\lambda| < \rho(C)$.
  \end{theorem}
\end{color}

\newpage
    \begin{color}{red}\LARGE\boldmath
      \section{\sc\bf The eigenvalues of the affine Coxeter transformation are roots of
unity}
      \end{color}
    \vspace{3mm}

The Coxeter transformation corresponding to the extended Dynkin
diagram is called the \begin{color}{blue}{\bf affine Coxeter
transformation}.\end{color} \vspace{2mm} \\

\begin{color}{blue}\LARGE\boldmath
\begin{theorem}
 The eigenvalues of the affine Coxeter transformation are roots of unity.
\end{theorem}
\end{color}

%%\DL{sm takzhe upr 26 na str 55 Koli B.(1968) (v dannom slu -- eshche
%% bolee rannij rez. nekogo Zho.-Py. Serra) }
  Subbotin-Stekolshchik \cite{SuSt79}, \cite{St82a}.
 The same theorem for the case of the Dynkin diagrams is due to Coxeter, \cite{Cox51},
 \cite{Cox49}.\\

\large The citation from \cite{Cox51}:
 \lq\lq
 Having computed the \begin{color}{blue}\boldmath$m$\end{color}'s
 several years earlier \cite{Cox49}, I recognized
them  in the Poincar\'{e} polynomials  while listening to
Chevalley's address at the International Congress in 1950. I am
grateful to A.~J.~Coleman for drawing my attention to the relevant
work of Racah, which helps to explain the ``coincidence''; also, to
J.~S.~Frame for many helpful suggestions... \rq\rq \vspace{2mm}
\Large

%%\DL{nehorosho:
%% G -- eto gruppa,
%% a gruppa Weyl'y ee algebry Lie -- W}
In this case: eigenvalues are as follows:
\begin{color}{blue}\boldmath
\begin{equation*}
   \omega^{m_1}, \quad \omega^{m_2}, \quad \dots, \quad \omega^{m_n},
\end{equation*}
\end{color} \\
where $\omega = e^{2{\pi}i/h}$,
\begin{color}{blue}\boldmath$h$\end{color} is the Coxeter number,
\begin{color}{blue}\boldmath$m_i$\end{color} are exponents of eigenvalues,
\begin{color}{blue}\boldmath$m_i + 1$\end{color} are the degrees of
 homogeneous basic elements of \begin{color}{blue}\boldmath$R^G$\end{color} is the
 \begin{color}{blue}{\bf algebra of invariants}\end{color} of the Weyl
 group $G$.

%%\DL{net slov:}
Let \begin{color}{blue}\boldmath$P(\mathcal{L}, t)$\end{color} be
 the
\begin{color}{blue}{\bf Poincar\'{e} series}\end{color} of the
corresponding
\begin{color}{blue}{\bf Lie group
\boldmath$\mathcal{L}$}\end{color}.  Then

\begin{color}{blue}\boldmath
\begin{equation*}
   P(\mathcal{L}, t)  = (1 + t^{2m_1 + 1})(1 + t^{2m_2 + 1})\dots(1 + t^{2m_n + 1}).
\end{equation*}
\end{color} \\
  (Hopf's theorem)
  \cite{CE48}, \cite{Col58}, \cite{Sol63}.

\newpage
~\\
    \begin{color}{red}\LARGE\boldmath
      \section{\sc\bf Splitting along the edge formula (Subbotin-Sumin)}
      \end{color}
    \vspace{3mm}

An edge  $l$ is said to be
 \begin{color}{blue}{\bf splitting}\end{color}
  if by deleting it we split the graph
$\varGamma$ into two graphs $\varGamma_1$ and $\varGamma_2$.
\vspace{2mm} \\

\begin{figure}[h]
\centering
\includegraphics{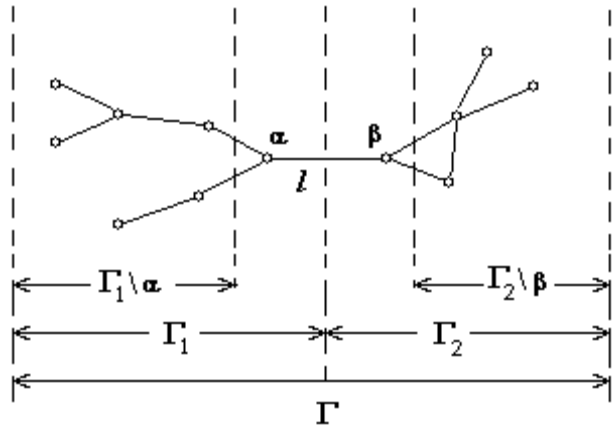}
\caption{\hspace{3mm} A split graph $\varGamma$}
%% splitted graph => split graph
%%%%%% The label must come after caption
\label{splitted_diagram}
\end{figure}

\begin{proposition}
  \label{split_edge}
  For a given graph $\varGamma$ with a splitting edge $l$, we have
  \vspace{2mm} \\
\begin{color}{blue}\Large\boldmath
 \begin{equation}
  \label{split_edge_eq}
   \mathcal{X}(\varGamma, \lambda) =
       \mathcal{X}(\varGamma_1, \lambda)\mathcal{X}(\varGamma_2, \lambda) -
       \lambda\mathcal{X}(\varGamma_1\backslash\alpha, \lambda)
              \mathcal{X}(\varGamma_2\backslash\beta, \lambda),
\end{equation}
\end{color}   \vspace{2mm} \\
where $\alpha$ and $\beta$ are the endpoints of the deleted edge
$l$.
\end{proposition}

Subbotin-Sumin \cite{SuSum82}. This is the simply-laced case.

\newpage
~\\
    \begin{color}{red}\LARGE\boldmath
      \section{\sc\bf Splitting along the edge formula (multiply-laced case)}
      \end{color}
    \vspace{3mm}

\begin{proposition}
  \label{split_edge_2}
  For a given graph $\varGamma$ with a splitting weighted edge $l$
  corresponding to roots of different lengths, we have

\begin{color}{blue}\Large\boldmath
\begin{equation*}
  \label{split_edge_eq_2}
   \mathcal{X}(\varGamma, \lambda) =
       \mathcal{X}(\varGamma_1, \lambda)\mathcal{X}(\varGamma_2, \lambda) -
       \rho\lambda\mathcal{X}(\varGamma_1\backslash\alpha, \lambda)
              \mathcal{X}(\varGamma_2\backslash\beta, \lambda),
\end{equation*}
\end{color} \\
  where $\alpha$ and $\beta$ are the endpoints
  of the deleted edge $l$, and $\rho$ is the following factor:
 \begin{equation*}
     \rho = k_{\alpha\beta} k_{\beta\alpha},
 \end{equation*}
   where $k_{ij}$ is an element
   of the Cartan matrix, see above examples $\widetilde{F}_{41}$,
   $\widetilde{F}_{42}$.
\end{proposition}\vspace{5mm}

\begin{corollary}
  \label{split_edge_corol_2}
  Let $\varGamma_2$ (in Proposition \ref{split_edge_2})
  be a component containing a single point. Then, the
  following formula holds
\begin{color}{blue}\Large\boldmath
\begin{equation*}
 \label{split_edge_eq_corol_2}
   \mathcal{X}(\varGamma, \lambda) =
       -(\lambda + 1)\mathcal{X}(\varGamma_1, \lambda) -
       \rho\lambda\mathcal{X}(\varGamma_1\backslash\alpha, \lambda),
\end{equation*}
\end{color}
\end{corollary}

\newpage
~\\
    \begin{color}{red}\LARGE\boldmath
      \section{\sc\bf Gluing formulas}
      \end{color}
    \vspace{3mm}

\begin{proposition}
 \label{gluing}
   Let $*_n$ be a star with $n$ rays coming from the vertex.
 Let $\varGamma(n)$ be the graph obtained from  $*_n$ by gluing $n$
 copies of the graph $\varGamma$ to the endpoints of its rays . Then

\begin{color}{blue}\Large\boldmath
\begin{equation*}
  \begin{split}
   & \mathcal{X}(\varGamma(n), \lambda) =
       \mathcal{X}(\varGamma, \lambda)^{n-1}\varphi_{n-1}(\lambda),
       \text { where } \\
   & \\
   & \varphi_{n}(\lambda) =
       \mathcal{X}(\varGamma + \beta, \lambda)
        - n\lambda\mathcal{X}(\varGamma\backslash\alpha, \lambda).
  \end{split}
\end{equation*}
\end{color}

\end{proposition}

 Subbotin-Sumin \cite{SuSum82}. (See, also \S\ref{sect_blocks}).

\begin{figure}[h]
\centering
\includegraphics{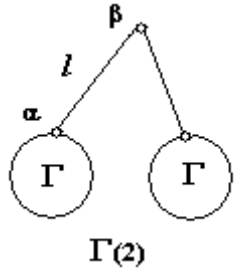}
\caption[\hspace{3mm} Splitting along the edge $l$ of the graph
$\Gamma(2)$]{\hspace{3mm} Splitting along the edge $l$ of the graph
$\Gamma(2)$. \\  Here, the graph $\Gamma(2)$ is obtained by gluing
two copies of the graph $\Gamma$.}
 \label{split_along_l}
\end{figure}

\begin{proposition}
  \label{merge_2G}
  If the spectrum of the Coxeter transformations
  for graphs $\varGamma_1$ and $\varGamma_2$ contains an eigenvalue $\lambda$,
  then this eigenvalue is also the eigenvalue of the Coxeter
  transformation for the graph $\varGamma$
  obtained by gluing as described in Proposition \ref{gluing}.
\end{proposition}

This proposition follows from the following formula:

\begin{color}{blue}\Large\boldmath
\begin{equation*}
  \begin{split}
    & \mathcal{X}(\varGamma_1 + \beta + \varGamma_2, \lambda) = \\
    & \mathcal{X}(\varGamma_1, \lambda)\mathcal{X}(\varGamma_2 + \beta, \lambda) -
    \lambda\mathcal{X}(\varGamma\backslash\alpha, \lambda)
    \mathcal{X}(\varGamma_2, \lambda).
  \end{split}
\end{equation*}
\end{color}

\newpage
~\\
    \begin{color}{red}\LARGE\boldmath
      \section{\sc\bf The Dynkin diagram ${A}_n$, the Frame formula}
      \end{color}
    \vspace{3mm}

\begin{color}{blue}\Large\boldmath
\begin{equation*}
 \label{recur_chi}
  \begin{split}
    \mathcal{X}({A}_1) = & -(\lambda + 1),  \\
    \mathcal{X}({A}_2) = & \lambda^2 + \lambda + 1,  \\
    \mathcal{X}({A}_3) = & -(\lambda^3 + \lambda^2 + \lambda + 1),  \\
    \mathcal{X}({A}_4) = & \lambda^4 + \lambda^3 + \lambda^2 + \lambda + 1,  \\
    \dots  \\
    \mathcal{X}({A}_n) = & -(\lambda + 1)\mathcal{X}({A}_{n-1}) -
           \lambda\mathcal{X}({A}_{n-2}),
\hspace{5mm} n > 2.
  \end{split}
\end{equation*}
\end{color}

J.~S.~Frame in \cite[p.784]{Fr51} obtained that

\begin{color}{blue}\Large\boldmath
\begin{equation*}
   \label{frame}
   \mathcal{X}({A}_{m+n}) =
     \mathcal{X}({A}_{m})\mathcal{X}({A}_{n})
   - \lambda\mathcal{X}({A}_{m-1})\mathcal{X}({A}_{n-1}),
\end{equation*}
\end{color}
 which easily follows from eq. (\ref{split_edge_eq}).

\newpage
    \begin{color}{red}\LARGE\boldmath
      \section{\sc\bf The spectral radius and Lehmer's number (McMullen)}
      \end{color}
    \vspace{1mm}

\begin{color}{blue}\Large\boldmath
\begin{theorem}
   Either
   $\rho({\bf C}) = 1$, or $\rho({\bf C}) \geq \lambda_{Lehmer} \approx 1.176281...$
   The spectral radius
$\rho({\bf C})$ of the Coxeter transformation for all graphs with
indefinite Tits form attains its minimum when the diagram is
${E}_{10}$.
\end{theorem}
\end{color}
 (McMullen, \cite{McM02}). \\

\begin{color}{blue}{\bf Lehmer's number}\end{color} is a root
$\mathcal{X}(C)$ for the diagram $E_{10}$.

\begin{color}{blue}\Large\boldmath
\begin{equation*}
   \mathcal{X}(C) = x^{10} + x^9 - x^7 - x^6 - x^5 - x^4 - x^3 + x + 1,
\end{equation*}
\end{color}
 \begin{figure}[h] \centering
\includegraphics{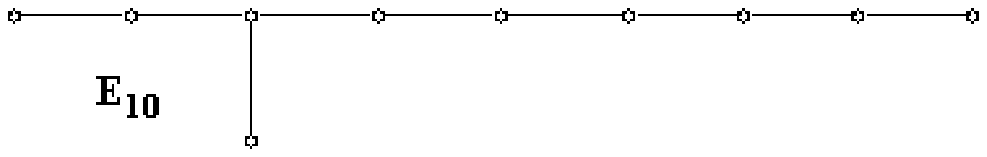}
%%\caption{\hspace{3mm} The spectral radius attains its minimum on the
%%graph ${E}_{10}$}
%%%%%% The label must come after caption
\label{E10_diagram}
\end{figure}

  Let $p(x)$ be a monic integer polynomial, and define its
\begin{color}{blue}{\bf Mahler measure}\end{color} to be

\begin{color}{blue}\Large\boldmath
\begin{equation*}
\Vert p(x)\Vert  = \prod\limits_{\beta} |\beta| ,
\end{equation*}
 \end{color} \\
where $\beta$ runs over all (complex) roots of $p(x)$ outside the
unit circle.

 %% \DL{tut nado vstavit' ego
 %% imya, uzh ne $\alpha$ li?}
 %%
 %% \DL{nichego ne ponimayu: $\Vert p(x)\Vert$
 %% opredelena dlya mnogochlenov, a $\alpha$ -- chislo}
 %% Kak dlya \alpha isp. $\Vert p(x)\Vert$ Dat' otvet in v.2
 In 1933, Lehmer \cite{Leh33} asks whether,
 for each $\varepsilon \ge 1$, there exists an algebraic
integer $\alpha$ such that

\begin{color}{blue}\Large\boldmath
\begin{equation}
 \label{numb_lehmer}
    1 < \Vert \alpha \Vert < 1 + \varepsilon.
\end{equation}
\end{color} \\
 %%The polynomial with minimal root $\alpha$ (in the
%%sense of (\ref{numb_lehmer})) Lehmer could find in \cite{Leh33} is
 %%$E10$.
%%\DL{vidimo, ty imel v vidu}
 In \cite{Leh33}, Lehmer established that the polynomial
 with minimal root $\alpha$ (in the sense of (\ref{numb_lehmer})) is
 $E10$. For details, see \cite{Hir02}.

\newpage
    \begin{color}{red}\LARGE\boldmath
      \section{\sc\bf The spectral radius of diagrams $T_{2,3,n}$
        and the Pisot number (Zhang)}
      \end{color}
    \vspace{1mm}

The following diagrams belong to the class
 $T_{2,3,n}$:   $D_5$ ($n=2$), $E_6$ ($n=3$), $E_7$ ($n=4$),
$E_8$ ($n=5$),
 $\widetilde{E}_8$ ($n = 6$), $E_{10}$ ($n = 7$). \vspace{2mm} \\

 \begin{proposition}
 \label{polyn_T_1}
The characteristic polynomials of Coxeter transformations for the
diagrams $T_{2,3,n}$ are as follows:
\begin{color}{blue}\LARGE\boldmath
\begin{equation*}
 \label{E_series}
  \mathcal{X}(T_{2,3,n-3}) =
     \lambda^n + \lambda^{n-1} -
     \sum\limits_{i=3}^{n-3}\lambda^{i} + \lambda + 1.
\end{equation*}
\end{color}

The spectral radius $\rho(T_{2,3,n-3})$ converges to the maximal
 root $\rho_{max}$
 of the equation
 \begin{color}{blue}\LARGE\boldmath
\begin{equation*}
 \label{zhang_eq}
  \lambda^3 - \lambda - 1 = 0,
\end{equation*}
 \end{color}

and
 \begin{color}{blue}\LARGE\boldmath
\begin{equation*}
 \label{zhang_number}
   \rho_{max} =
   \sqrt[3]{\frac{1}{2} + \sqrt{\frac{23}{108}}} +
   \sqrt[3]{\frac{1}{2} - \sqrt{\frac{23}{108}}} \approx 1.324717...\ .
\end{equation*}
 \end{color}
\end{proposition}

 The fact that \begin{color}{blue}\Large\boldmath$\rho(T_{2,3,n}) \rightarrow  \rho_{max}$\end{color}
 as \begin{color}{blue}\Large\boldmath$n \rightarrow \infty$\end{color}
  was obtained by Zhang \cite{Zh89} and used in the study of regular components
  of an Auslander-Reiten quiver. The number
   \begin{color}{blue}\Large\boldmath$\rho_{max}$\end{color}
  coincides with \begin{color}{blue}\Large{\bf Pisot
  number}\end{color}.

Recall that an algebraic integer $\lambda > 1$ is said to be a Pisot
number if all its conjugates (other then $\lambda$ itself) satisfy
$|\lambda^{'}| < 1$.

The smallest Pisot number is a root of $\lambda^3 - \lambda - 1 = 0$
:
\begin{color}{blue}\LARGE\boldmath
\begin{equation*}
  \label{pisot_zhang}
          \lambda_{Pisot} \approx 1.324717...
\end{equation*}
\end{color}

\newpage
~\\ \vspace{2mm}
    \begin{color}{red}\LARGE\boldmath
      \section{\sc\bf The spectral radii of the diagrams $T_{3,3,n}$}
      \end{color}
    \vspace{1mm}

Recall that the diagrams $E_6$ ($n=2$) and $\widetilde{E}_6$ ($n=3$)
belong to the class
 $T_{3,3,n}$. \vspace{2mm} \\

\begin{proposition}
 \label{polyn_T_2}
The characteristic polynomials of Coxeter transformations for %%%${\bf E}_n$-series,
the  diagrams $T_{3,3,n}$ with $n \geq 3$ are as follows:

\begin{color}{blue}\Large\boldmath
\begin{equation*}
 \label{E6_series}
  \mathcal{X}(T_{3,3,n}) =
     \lambda^{n+4} + \lambda^{n+3}
     -2\lambda^{n+1}
     -3\sum\limits_{i=4}^{n}\lambda^{i}
     -2\lambda^3 + \lambda + 1,
\end{equation*}
\end{color}

The spectral radius $\rho(T_{3,3,n})$ converges to the maximal
 root $\rho_{max}$
 of the equation

\begin{color}{blue}\Large\boldmath
\begin{equation*}
 \label{spectr_rad_3}
  \lambda^2 - \lambda - 1 = 0,
\end{equation*}
\end{color}

and

\begin{color}{blue}\Large\boldmath
\begin{equation*}
 \label{spectr_rad_4}
   \rho_{max} =
   \frac{\sqrt{5} + 1}{2} \approx 1.618034...\ \text{ (the Golden mean) }.
\end{equation*}
\end{color}

\end{proposition}

\newpage
~\\ \vspace{2mm}
    \begin{color}{red}\LARGE\boldmath
      \section{\sc\bf The spectral radii of the diagrams $T_{2,4,n}$ (Lakatos) }
      \end{color}
    \vspace{1mm}

Recall that the diagrams ${D}_6 (n = 2)$, ${E}_7 (n = 3)$, and
  $\widetilde{E}_7 (n = 4)$ belong to the class
 $T_{2,4,n}$. \vspace{2mm} \\

 \begin{proposition}
 \label{polyn_T_3}
The characteristic polynomials of Coxeter transformations for
diagrams $T_{2,4,n}$, where $n \geq 3$,
  are as follows:
\begin{color}{blue}\Large\boldmath
\begin{equation*}
 \label{E7_series}
  \mathcal{X}(T_{2,4,n}) =
     \lambda^{n+4} + \lambda^{n+3}
     -\lambda^{n+1}
     -2\sum\limits_{i=4}^{n}\lambda^{i}
     -\lambda^3 + \lambda + 1,
\end{equation*}
 \end{color}

The spectral radius $\rho(T_{2,4,n})$ converges to the maximal
 root $\rho_{max}$ of the equation

\begin{color}{blue}\Large\boldmath
\begin{equation*}
 \label{spectr_rad_1}
  \lambda^3 - \lambda^2 - 1 = 0,
\end{equation*}
 \end{color}

and

\begin{color}{blue}\Large\boldmath
\begin{equation*}
   \rho_{max} =
   \frac{1}{3} + \sqrt[3]{\frac{58}{108} + \sqrt{\frac{31}{108}}} +
   \sqrt[3]{\frac{58}{108} - \sqrt{\frac{31}{108}}} \approx 1.465571...\ .
\end{equation*}
 \end{color}

\end{proposition}

 Lakatos \cite{Lak99} obtained results on the convergence of the spectral radii
$\rho_{max}$ similar to propositions regarding $\rho(T_{2,3,n})$,
 $\rho(T_{3,3,n})$, $\rho(T_{2,4,n})$.

\newpage
~\\ \vspace{2mm}
    \begin{color}{red}\LARGE\boldmath
      \section{\sc\bf The binary polyhedral groups }
      \end{color}
    \vspace{1mm}

We consider the double covering

\begin{color}{blue}\LARGE\boldmath
\begin{equation*}
   \pi : SU(2) \longrightarrow SO(3,\mathbb{R}).
\end{equation*}
\end{color}

If $G$ is a finite subgroup of $SO(3,\mathbb{R})$, we see that the
preimage $\pi^{-1}(G)$ is a finite subgroup of $SU(2)$ and
$|\pi^{-1}(G)| = 2|G|$. The finite subgroups of $SO(3,\mathbb{R})$
are called \begin{color}{blue}{\bf polyhedral groups}\end{color},
see \text{Table \ref{rotation_pol_1}}. The finite subgroups of
$SU(2)$ are naturally called \begin{color}{blue}{\bf binary
polyhedral groups}\end{color}, see Table \ref{binary_pol_1}.

\begin{table} [h]
  \centering
  \vspace{2mm}
  \caption{\hspace{3mm}The polyhedral groups in $\mathbb{R}^3$}
  \renewcommand{\arraystretch}{1.5}
 \begin{color}{blue}\Large\boldmath
 \begin{tabular} {||c|c|c|c|c||}
 \hline \hline
   Polyhedron & Orders of symmetries & Rotation group & Group order \\
 \hline \hline
     Pyramid   &  $-$  & cyclic & $n$ \\
 \hline
     Dihedron  &  $n$  $2$ $2$  & dihedral & $2n$ \\
 \hline
     Tetrahedron &  $3$  $2$ $3$  & $\mathcal{A}_4$ & $12$ \\
 \hline
     Cube    &  $4$  $2$ $3$  & $\mathcal{S}_4$ & $24$ \\
     Octahedron  &  $3$  $2$ $4$  & $\mathcal{S}_4$ & $24$ \\
 \hline
     Dodecahedron &  $5$  $2$ $3$  & $\mathcal{A}_5$ & $60$ \\
     Icosahedron &  $3$  $2$ $5$  & $\mathcal{A}_5$ & $60$  \\
 \hline \hline
 \end{tabular}
 \end{color}
\\ \vspace{2mm}
 Here, $\mathcal{S}_m$ (resp. $\mathcal{A}_m$) denotes the
 symmetric, (resp. alternating) \\
 group of all (resp. of all even)
 permutations of $m$ letters.
 \label{rotation_pol_1}
\end{table}

\newpage
~\\ \vspace{2mm}
    \begin{color}{red}\LARGE\boldmath
      \section{\sc\bf The binary polyhedral groups (2) }
      \end{color}
    \vspace{1mm}

\begin{table} [h]
  \centering
  \vspace{2mm}
  \caption{\hspace{3mm}The finite subgroups of $SU(2)$}
  \renewcommand{\arraystretch}{1.5}
 \begin{color}{blue}\Large\boldmath
 \begin{tabular} {||c|c|c|c||}
 \hline \hline
   $\langle l, m, n \rangle$ & Order & Notation & Well-known name \\
 \hline \hline
      $-$            &  $n$ & $\mathbb{Z}/n\mathbb{Z}$ & cyclic
group \\
 \hline
     $\langle 2, 2, n \rangle$ &    $4n$     & $\mathcal{D}_n$ & binary
dihedral group \\
 \hline
     $\langle 2, 3, 3 \rangle$ &    $24$  & $\mathcal{T}$ & binary tetrahedral
group \\
 \hline
     $\langle 2, 3, 4 \rangle$ &    $48$   & $\mathcal{O}$ & binary
octahedral group \\
 \hline
     $\langle 2, 3, 5 \rangle$ &    $120$  & $\mathcal{J}$ & binary
icosahedral group \\
 \hline \hline
 \end{tabular}
 \end{color}
  \label{binary_pol_1}
\end{table}

The \underline{binary polyhedral group} is generated by three
generators $R$, $S$, and $T$ subject to the relations
 \begin{color}{blue}\LARGE\boldmath
\begin{equation*}
  \label{natural_gen}
   R^p = S^q = T^r = RST = -1.
\end{equation*}
\end{color}

Denote  this group by $\langle p, q, r \rangle$. The order of the
group $\langle p, q, r \rangle$ is
\begin{color}{blue}\LARGE\boldmath
\begin{displaymath}
  \displaystyle\frac{4}
   {\displaystyle\frac{1}{p} + \frac{1}{q} + \frac{1}{r} - 1}.
\end{displaymath}
\end{color}

\newpage
~\\ \vspace{2mm}
    \begin{color}{red}\LARGE\boldmath
      \section{\sc\bf The binary polyhedral groups, the algebra of invariants (F.~Klein) }
      \end{color}
    \vspace{1mm}

\begin{theorem}
The algebra of invariants $\mathbb{C}[z_1, z_2]^G$ is generated by 3
indeterminates $x, y, z,$ subject to one %essential
relation

\begin{color}{blue}\LARGE\boldmath
\begin{equation}
 \label{R_curve}
  R(x,y,z) = 0,
\end{equation}
\end{color} \\
where $R(x,y,z)$ is defined in Table \ref{klein_rel}. In other
words, the algebra of invariants $\mathbb{C}[z_1, z_2]^G$ coincides
with the coordinate algebra of the curve defined by Eq.
(\ref{R_curve}), i.e.,

\begin{color}{blue}\LARGE\boldmath
\begin{equation}
 \label{quotient_G}
 \mathbb{C}[z_1, z_2]^G \simeq \mathbb{C}[x, y, z]/(R(x,y,z)).
\end{equation}
\end{color}

\end{theorem}

  F.~Klein, 1884, \cite{Kl1884}.

\begin{table} [h]
  \centering
  \vspace{2mm}
  \caption{\hspace{3mm}The relations $R(x, y, z)$ describing the algebra of invariants
      $\mathbb{C}[z_1, z_2]^G$}
  \renewcommand{\arraystretch}{1.5}
 \begin{color}{blue}\Large\boldmath
 \begin{tabular} {||c|c|c||}
 \hline \hline
    Finite subgroup of $SU(2)$ & Relation $R(x, y, z)$ & Dynkin diagram \\
 \hline \hline
      $\mathbb{Z}/n\mathbb{Z}$ & $x^n + yz$ & $A_{n-1}$ \\
 \hline
      $\mathcal{D}_n$ & $x^{n+1} + xy^2 + z^2$ & $D_{n+2}$ \\
 \hline
      $\mathcal{T}$ & $x^4 + y^3 + z^2$ & $E_6$ \\
 \hline
      $\mathcal{O}$ & $x^3{y} + y^3 + z^2$ & $E_7$ \\
 \hline
      $\mathcal{J}$ & $x^5 + y^3 + z^2$ & $E_8$ \\
 \hline \hline
 \end{tabular}
 \end{color}
  \label{klein_rel}
\end{table}

\newpage
~\\ \vspace{2mm}
    \begin{color}{red}\LARGE\boldmath
      \section{\sc\bf The binary polyhedral groups, Kleinian singularities }
      \end{color}
    \vspace{1mm}

%% \DL{tak chto takoe (\ref{quotient_G}): algebra (kak
%% ya by skazal) ili variety, na kotorom eta algebra -- algebra f-ij?}
%% Otvetit' in v.2
The quotient algebra (\ref{quotient_G}) has no singularity except at
the origin $O \in \mathbb{C}^3$. The quotient variety (or, orbit
space)
 $X = \mathbb{C}^2/G$ is isomorphic to (\ref{quotient_G}) (see,
 \cite{Hob02}).

The quotient variety  $X$
 is called a \begin{color}{blue}{\bf Kleinian singularity}\end{color}
  also known as a \begin{color}{blue}{\bf Du Val
  singularity}\end{color}.

\newpage
~\\ \vspace{2mm}
    \begin{color}{red}\LARGE\boldmath
      \section{\sc\bf The binary polyhedral groups, algebras of invariants. An example }
      \end{color}
    \vspace{1mm}

Consider the cyclic group $G = \mathbb{Z}/r\mathbb{Z}$ of order $r$.
The group $G$ acts on $\mathbb{C}[z_1, z_2]$ as follows:

 \begin{color}{blue}\LARGE\boldmath
    \begin{equation*}
       (z_1, z_2)  \mapsto (\varepsilon{z_1},
       \varepsilon^{r-1}z_2),
    \end{equation*}
  \end{color} \\
   where $\varepsilon = e^{2\pi{i}/r}$, and the polynomials

  \begin{color}{blue}\LARGE\boldmath
    \begin{equation*}
        x = z_1z_2, \hspace{5mm} y = -z_1^r, \hspace{5mm} z = z_2^r
    \end{equation*}
  \end{color} \\
  are invariant polynomials in $\mathbb{C}[x, y, z]$ which
  satisfy the following relation

  \begin{color}{blue}\LARGE\boldmath
    \begin{equation*}
       x^r + yz = 0,
    \end{equation*}
  \end{color}

  We have
  \begin{color}{blue}\LARGE\boldmath
    \begin{equation*}
       k[V]^G = \mathbb{C}[z_1z_2, z_1^r, z_2^r]
                \simeq \mathbb{C}[x, y, z]/(x^r + yz) .
    \end{equation*}
  \end{color}

\newpage
~\\ \vspace{2mm}
    \begin{color}{red}\LARGE\boldmath
      \section{\sc\bf The binary polyhedral groups,\\ connection with Dynkin diagrams \\ (Du Val's phenomenon) }
      \end{color}
    \vspace{1mm}

Du Val obtained the following description of the minimal resolution
\begin{color}{blue}\LARGE\boldmath
    \begin{equation*}
  \pi: \tilde{X} \longrightarrow X
    \end{equation*}
  \end{color}
of a Kleinian singularity $X = \mathbb{C}^2/G$, \cite{DuVal34}

The {\it exceptional divisor} (the preimage of the singular point
$O$) is a finite union of complex projective lines:

\begin{color}{blue}\LARGE\boldmath
    \begin{equation*}
   \pi^{-1}(O) = L_1 \cup \dots \cup L_n,
   \hspace{3mm} L_i \simeq \mathbb{CP}^1 \text{ for } i =
   1,\dots,n.
    \end{equation*}
  \end{color}
For $i \neq j$, the intersection $L_i \cap L_j$ is empty or consists
of exactly one point.

To each complex projective line $L_i$ (which can be identified with
the sphere $S^2 \subset \mathbb{R}^3$) we assign a vertex $i$, and
two vertices are connected by an edge if the corresponding
projective lines intersect. The corresponding diagrams are Dynkin
diagrams, see \text{Table \ref{klein_rel}}.

\newpage
~\\ \vspace{2mm}
    \begin{color}{red}\LARGE\boldmath
      \section{\sc\bf The binary polyhedral groups, Du Val's phenomenon for binary dihedral group }
      \end{color}
    \vspace{1mm}

\begin{figure}[h]
\centering
\includegraphics{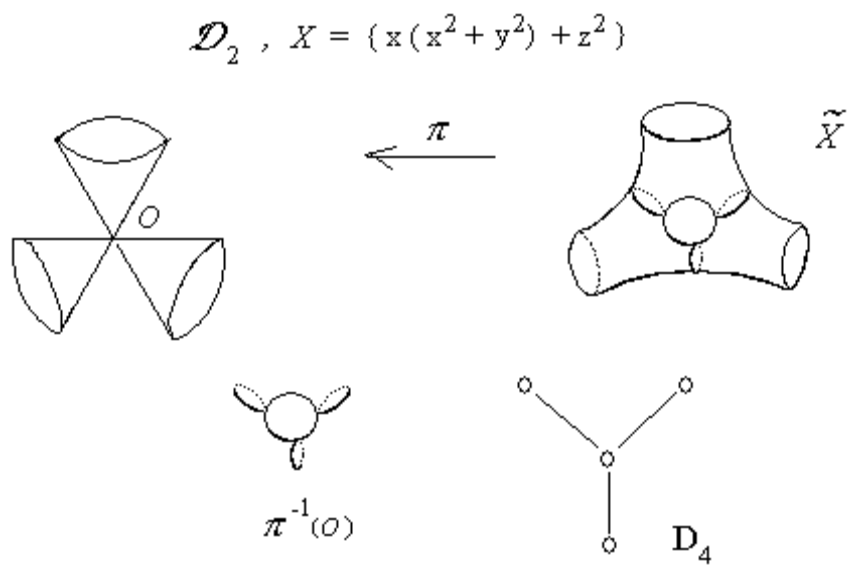}
\caption{\hspace{3mm} The minimal resolution $\pi^{-1}: \tilde{X}
\longrightarrow X$
  for $X = \mathcal{D}_2$}
%%%%%% The label must come after caption
\label{D2_singularity}
\end{figure}

%% \DL{rather graphic?} - soglaseb; maybe: visual picture; obvious pictuare ?
In the case of the binary dihedral group $\mathcal{D}_2$, the real
resolution of the real variety
$$
   \mathbb{C}^3/R(x,y,z) \cap \mathbb{R}^3
$$
gives a rather graphic picture of the complex situation, the minimal
resolution $\pi^{-1}: \tilde{X} \longrightarrow X$ for $X =
\mathcal{D}_2$ is depicted on Fig.~\ref{D2_singularity}. Here
$\pi^{-1}(O)$ consists of four circles, the corresponding diagram is
the Dynkin diagram $D_4$.

\newpage
~\\ \vspace{2mm}
    \begin{color}{red}\LARGE\boldmath
      \section{\sc\bf The McKay correspondence }
      \end{color}
    \vspace{1mm}

Let $G$ be a finite subgroup of $SU(2)$. Let $\{\rho_0,\rho_1,\dots,
\rho_n\}$ be the set of all distinct irreducible finite dimensional
complex representations of $G$, of which $\rho_0$ is the trivial
one. Let $\rho: G \longrightarrow SU(2)$ be a faithful
 representation, then, for each group $G$, we define
 a matrix $A(G) = (a_{ij})$, by decomposing the tensor products:

\begin{color}{blue}\Large\boldmath
\begin{equation}
 \label{main_McKay}
 \rho \otimes \rho_j = \bigoplus\limits_{k=0}^r a_{jk}\rho_k,
  \hspace{7mm} j = 0,1,...,r,
\end{equation}
\end{color}

where $a_{jk}$ is the multiplicity of $\rho_k$ in $\rho \otimes
\rho_j$. McKay observed that \\

\begin{color}{blue}\Large\boldmath
\begin{center}
{\it The matrix $2I - A(G)$ is the Cartan matrix of the extended
  Dynkin diagram $\tilde\Gamma(G)$ associated to $G$.
  There is a one-to-one correspondence between finite subgroups
  of $SU(2)$ and simply-laced extended Dynkin diagrams.
 }\end{center}
\end{color} \vspace{3mm}

For the multiply-laced case, the McKay correspondence was extended
by D.~Happel, U.~Preiser, and C.~M.~Ringel, \cite{HPR80} and by
P.~Slodowy, \cite{Sl80}. We consider P.~Slodowy's approach.

The systematic proof of the McKay correspondence based on the study
of affine Coxeter transformations was given by R. Steinberg,
\cite{Stb85}.

\newpage
~\\ \vspace{2mm}
    \begin{color}{red}\LARGE\boldmath
      \section{\sc\bf The Slodowy correspondence }
      \end{color}
    \vspace{1mm}

Slodowy's approach is based on the consideration of
\begin{color}{blue}{\bf restricted representations}\end{color} and
\begin{color}{blue}{\bf induced representations}\end{color} instead of an
original representation. Let $\rho: G \longrightarrow GL(V)$ be a
representation of a group $G$. We denote the \underline{restricted
representation} of $\rho$ to a subgroup $H \subset G$ by
$\rho\downarrow^G_H$, or, briefly, $\rho^\downarrow$ for fixed $G$
and $H$. Let $\tau: H \longrightarrow GL(V)$ be a representation of
a subgroup $H$. We denote by $\tau\uparrow^G_H$ the representation
\underline{induced} by $\tau$ to a representation of the group $G$
containing $H$; we briefly write $\tau^\uparrow$ for fixed $G$ and
$H$.

Let us consider pairs of groups $H \triangleleft G$, where $H$ and
$G$ are binary polyhedral groups from Table \ref{group_pairs}.

\begin{table} [h]
  \centering
  \vspace{2mm}
  \caption{\hspace{3mm}The pairs $H \triangleleft G$ of binary polyhedral groups}
  \renewcommand{\arraystretch}{0.5} %% 1.5
 \begin{color}{blue}\boldmath
 \begin{tabular} {||c|c|c|c|c||}
 \hline \hline
   \quad Subgroup \quad &
   \quad Dynkin  \quad  &
   \quad Group \quad    &
   \quad Dynkin \quad   &
   \quad Index \quad    \cr
   \quad $H$  \quad     &
   \quad diagram \quad &
   \quad $G$  \quad &
   \quad diagram \quad &
   \quad $[G:H]$ \quad \cr
   \quad &
   \quad $\varGamma(H)$ &
   \quad &
   \quad $\varGamma(G)$ & \\
 \hline \hline
      & & & & \cr
      $\mathcal{D}_2$ & ${D}_4$ & $\mathcal{T}$ & ${E}_6$ & $3$ \cr
      & & & & \\
  \hline
      & & & & \cr
      $\mathcal{T}$  & ${E}_6$ & $\mathcal{O}$ & ${E}_7$ & $2$ \cr
      & & & & \\
  \hline
      & & & & \cr
      $\mathcal{D}_{n-1}$ & ${D}_{n+1}$ & $ \mathcal{D}_{2(n-1)}$ &
      ${D}_{2n}$ & $2$ \cr
      & & & & \\
 \hline
      & & & & \cr
      $\mathbb{Z}/2n\mathbb{Z}$ & ${A}_{2n-1}$ & $\mathcal{D}_n$ &
      ${D}_{n+2}$ & $2$ \cr
      & & & & \\
 \hline \hline
 \end{tabular}
 \end{color}
  \label{group_pairs}
\end{table}

Let us fix a pair $H \triangleleft G$ from Table \ref{group_pairs}.
We formulate now the essence of the Slodowy correspondence.

\newpage
~\\ \vspace{2mm}
    \begin{color}{red}\LARGE\boldmath
      \section{\sc\bf Induced representations; an example }
      \end{color}
    \vspace{1mm}
%% \DL{non-tirivial?}
Let $G$ be a finite group and $H$ any subgroup of $G$. Let
\begin{color}{blue}\boldmath$\tau$\end{color} be a representation of $H$ in the vector
space $V$.  The induced representation
\begin{color}{blue}\boldmath$\tau\uparrow^G_H$\end{color} of $G$
(or, \begin{color}{blue}\boldmath$\tau^\uparrow$\end{color}, or
\begin{color}{blue}\boldmath$\operatorname{Ind}_H^G\tau$\end{color})
in the space
\begin{color}{blue}\boldmath
  \begin{equation}
   \label{induced_1}
     W=\bigoplus_{x\in G/H} xV
  \end{equation}
  \end{color}
is defined as follows:
 \begin{color}{blue}\boldmath
 \begin{equation}
   \label{induced_2}
       g\cdot\sum_{x\in G/H} x v_x = \sum_{x\in G/H} gx v_x,
 \end{equation}
 \end{color} \\
 where $v_x \in V$ for each $x$. \vspace{2mm}\\

 \underline{Example.} Let $H$ be a cyclic group of order $3$,
 $H = \{1, a, a^2 \}$. Let $\omega := e^{{2\pi{i}/3}}$.
 There are $3$ irreducible representations of $H$, or $3$
 irreducible $\mathbb{C}H$-submodules of $\mathbb{C}H$:

\begin{color}{blue}\boldmath
\begin{equation*}
  \begin{split}
     & \tau_0 = \{ 1 + a + a^2 \} ; \quad a*z = z \\
     \\
     & \tau_1 = \{ 1 + \omega^2{a}  + \omega{a^2}\} ; \quad a*z = \omega{z} \\
     \\
     & \tau_2 = \{ 1 + \omega{a} + \omega^2{a^2} \} ; \quad a*z = \omega^2{z} \\
  \end{split}
\end{equation*}
 \end{color}
and
\begin{color}{blue}\boldmath
 \begin{equation*}
       \mathbb{C}H = \tau_0 \oplus \tau_1 \oplus \tau_2.
 \end{equation*}
 \end{color} \\

\newpage
~\\ \vspace{2mm}
    \begin{color}{red}\LARGE\boldmath
      \section{\sc\bf Induced representations; an example (2) }
      \end{color}
    \vspace{1mm}

 Let $G$ be the rotation group of the triangle

\begin{color}{blue}\boldmath
\begin{equation*}
  \{a, b \mid  a^3 = b^2 = 1, ab = ba^2 \},
\end{equation*}
\end{color}

The three irreducible \underline{right $\mathbb{C}G$-submodules} of
$\mathbb{C}G$ are as follows:
\begin{color}{blue}\boldmath
\begin{equation*}
  \begin{split}
     & U_1 = \{ 1 + a + a^2 + b + ab + a^2{b} \}, \\
     & \text{ corresponding representation: }
      \rho_1: a \rightarrow  1, b \rightarrow  1, \\
     \\
     & U_2 = \{ 1 + a + a^2 - b - ab - a^2{b} \}, \\
     & \text{ corresponding representation: }
      \rho_2: a \rightarrow  1, b \rightarrow  -1,   \\
     \\
     & U_3 = \{ 1 + \omega^2{a} + \omega{a}^2, \quad
                b + \omega{b}a + \omega^2{b}{a}^2 \}, \\
     & U_4 = \{ 1 + \omega{a} + \omega^2{a}^2, \quad
                b + \omega^2{b}a + \omega{b}{a}^2 \}, \\
     & \text{ corresponding representation: } \\
     & \rho_3 : a \rightarrow
     \left (
      \begin{array}{cc}
        \omega & 0 \\
        0 & \omega^2 \\
      \end{array}
      \right ),
      b \rightarrow
      \left (
      \begin{array}{cc}
        0 & 1 \\
        1 & 0 \\
      \end{array}
      \right ). \\
      \\
      & \quad \mathbb{C}G = U_1 \oplus  U_2 \oplus U_3 \oplus U_4;
     \quad U_3 \simeq U_4. \\
  \end{split}
\end{equation*}
\end{color}

\newpage
~\\ \vspace{2mm}
    \begin{color}{red}\LARGE\boldmath
      \section{\sc\bf Induced representations; an example (3) }
      \end{color}
    \vspace{1mm}

%% \DL{strannoe nachalo. I kto takaya H?}
 Then, $H \subset G$, elements \{$1$, $b$\} are \underline{two left cosets}
 of $G/H$, and by (\ref{induced_1}), (\ref{induced_2})
 the induced representations of $G$ are as follows:

\begin{color}{blue}\boldmath
\begin{equation*}
  \begin{split}
     & \tau_0^\uparrow = \{ 1 + a + a^2, \quad b + ab + a^2{b} \} = \rho_1 \oplus \rho_2, \\
     \\
     & \tau_1^\uparrow = \{1 + \omega^2{a}  + \omega{a^2},
            \quad b + \omega^2a{b} + \omega{a}^2{b} \},   \\
     \\
     & \tau_2^\uparrow = \{ 1 + \omega{a} + \omega^2{a^2},
        \quad b + \omega{a}{b} + \omega^2{a}^2{b}
         \},  \\
     \\
     &  \tau_1^\uparrow  \simeq  \tau_2^\uparrow \simeq \rho_3. \\
  \end{split}
\end{equation*}
\end{color}

  Here, \begin{color}{blue}\boldmath
     $b + ab + a^2{b} = b + ba^2 + ba$
   \end{color}, and, equivalently, the right cosets may be considered.

\newpage
~\\ \vspace{2mm}
    \begin{color}{red}\LARGE\boldmath
      \section{\sc\bf The trivial representation, the Frobenius reciprocity }
      \end{color}
    \vspace{1mm}

A \begin{color}{blue}{\bf trivial representation}\end{color}
    is a representation $(V, \rho)$ of a group $G$
    on which all elements of $G$ act as the identity mapping of $V$.
    The character of the trivial representation %takes the value of one
is equal to $1$  at any group element.

    \begin{color}{blue}{\bf The Frobenius reciprocity.}\end{color}
    For characters of restricted representation $\psi^\downarrow = \psi\downarrow^G_H$ and
    the induced representation $\chi^\uparrow = \chi\uparrow^G_H$,  the following relation
    holds:

\begin{color}{blue}\boldmath
\begin{equation}
  \label{frob_1}
    \langle\psi, \chi^\uparrow\rangle_G ~=~
    \langle\psi^\downarrow, \chi\rangle_H.
\end{equation}
\end{color}

   Let us apply (\ref{frob_1}) to the \underline{trivial representation} $\psi$ of
   $G$.  Let $\chi$ be a \underline{non-trivial irreducible representation} of
   $H$.  Since $\psi^\downarrow$ is a trivial representation of $H$,
   we have $\langle\psi^\downarrow, \chi\rangle_H = 0$, and

\begin{color}{blue}\boldmath
\begin{equation}
  \label{frob_2}
    \langle\psi, \chi^\uparrow\rangle_G = 0.
\end{equation}
\end{color}

We will use (\ref{frob_2}) in the proof of the generalized Ebeling
theorem, see \S\ref{Ebeling_th}.

\newpage
    \begin{color}{red}\LARGE\boldmath
      \section{\sc\bf Restricted representations, Clifford's theorem }
      \end{color}
    \vspace{1mm}

\large
 See, \cite[\S20]{JL01}.
 In this section, we suppose $H \triangleleft G$.

\begin{theorem}[Clifford]
  Let $\chi$ be an irreducible  character of $G$. Then

  (1) all the constituents of $\chi^\downarrow_H$ have the same
  degree %% \DL{i.e., dimension};

  (2) if $\psi_1, \dots, \psi_m$ are all the constituents of the
  $\chi^\downarrow_H$, then for a positive integer $e$, we have

\begin{color}{blue}\boldmath
  \begin{equation*}
     \chi^\downarrow_H = e(\psi_1 + ... + \psi_m).
  \end{equation*}
\end{color}
\end{theorem}

In the following corollary from Clifford's theorem, we assume that \\
$[G:H] = 2$ (resp. $3$). We are interested in these cases, see Table
\ref{group_pairs}.

  \begin{proposition}
    \label{prop_irred_ch}
    Let $\chi$ be an irreducible character of $G$. Then either

    (1) $\chi^\downarrow_H$ is irreducible, or

    (2) $\chi^\downarrow_H$ is the sum of $2$ (resp. $3$) distinct irreducible
    characters of $H$ of the same degree. In this case, we have
 \begin{color}{blue}\boldmath
  \begin{equation*}
     \chi^\downarrow_H = \psi_1 + \psi_2, \text{ resp. }
     \chi^\downarrow_H = \psi_1 + \psi_2 + \psi_3. \\
  \end{equation*}
 \end{color}
   If $\psi$ is an irreducible character of $G$ such that
   $\psi^\downarrow_H$ has $\psi_1$ or $\psi_2$ (resp., or $\psi_3$)
   as a constituent, then \begin{color}{blue}\boldmath$\psi = \chi$\end{color}.
 \end{proposition}
\Large

   Let $\tilde\pi$ be the \underline{trivial representation} of $G$,
   and let $\chi^\downarrow_H$ be of case (2) from Prop. \ref{prop_irred_ch},
   and  $\tilde\pi \neq \chi$.
   Then \begin{color}{blue}\boldmath$\pi := \tilde\pi^\downarrow$\end{color}
   is the trivial representation of $H$,
   and $\pi$ does not contain $\psi_i$ as a
   constituent, and

 \begin{color}{blue}\boldmath
  \begin{equation}
    \label{clifford_2}
   \langle \pi, \chi^\downarrow_H \rangle = 0.
  \end{equation}
 \end{color}

We will use (\ref{clifford_2}) in the proof of the generalized
Ebeling theorem, see \S\ref{Ebeling_th}.
 \begin{remark}
  \label{rem_v0}
 {\rm
  For case (1) from Prop. \ref{prop_irred_ch}, there exist
  non-trivial irreducible representation
  \begin{color}{blue}\boldmath$\tilde\pi \neq \chi$\end{color} of $G$, such
  that \begin{color}{blue}\boldmath$\pi = \chi^\downarrow_H$\end{color}.
  Then, two representations \begin{color}{blue}\boldmath$\pi$\end{color} and
  \begin{color}{blue}\boldmath$\chi^\downarrow_H$\end{color}
  are gluing on %% \DL{prover': on ili kak ino} the vertex of
  the corresponding folded diagram associated
  with the Slodowy correspondence, \S\ref{Slodowy_2},
  \S\ref{folded_diag}.
  }
 \end{remark}

\newpage
    \begin{color}{red}\LARGE\boldmath
      \section{\sc\bf The Slodowy correspondence (2) }
        \label{Slodowy_2}
      \end{color}
    \vspace{1mm}

 1) Let $\rho_i$, where $i = 1,\dots,n$, be all
irreducible representations of $G$; let $\rho^\downarrow_i$ be the
corresponding restricted representations of the subgroup $H$. Let
$\rho$ be a faithful representation of $H$, which may be considered
as the restriction of a fixed faithful representation $\rho_{f}$ of
$G$. Then the following decomposition formula makes sense
\\ \vspace{2mm}
\begin{color}{blue}\Large\boldmath
\begin{equation}
  \label{matrix_Slodowy}
     \rho \otimes \rho^\downarrow_i =
          \bigoplus\limits_{j} a_{ji} \rho^\downarrow_j
\end{equation}
\end{color} \vspace{1mm}\\
and uniquely determines an $n\times{n}$ matrix $\widetilde{A} =
(a_{ij})$ such that

\begin{color}{blue}\Large\boldmath
\begin{equation}
  \label{lab_slodowy_1}
   K = 2I - \widetilde{A},
\end{equation}
\end{color} \vspace{1mm}\\
where $K$ is the Cartan matrix of the corresponding
\underline{folded} extended
Dynkin diagram. \vspace{1mm}\\

2) Let $\tau_i$, where $i = 1,\dots,n$, be all irreducible
representations of the subgroup $H$, let $\tau^\uparrow_i$ be the
induced representations of the group $G$. Then the following
decomposition formula makes sense \vspace{1mm}\\
\begin{color}{blue}\Large\boldmath
\begin{equation}
  \label{matrix_Slodowy_1}
     \rho \otimes \tau^\uparrow_i =
          \bigoplus a_{ij} \tau^\uparrow_j,
\end{equation}
\end{color} \vspace{1mm}\\
i.e., the decomposition of the induced representation is described
by the matrix $A^\vee = A^t$ which satisfies the relation

\begin{color}{blue}\Large\boldmath
\begin{equation}
    \label{lab_slodowy_2}
   K^\vee = 2I - \widetilde{A}^\vee,
\end{equation}
\end{color} \vspace{1mm}\\
where $K^\vee$ is the Cartan matrix of the dual \underline{folded}
 extended Dynkin diagram.

\newpage
~\\ \vspace{2mm}
    \begin{color}{red}\LARGE\boldmath
      \section{\sc\bf The Slodowy correspondence, folded diagrams }
        \label{folded_diag}
      \end{color}
    \vspace{1mm}

\begin{figure}[h]
\centering
\includegraphics{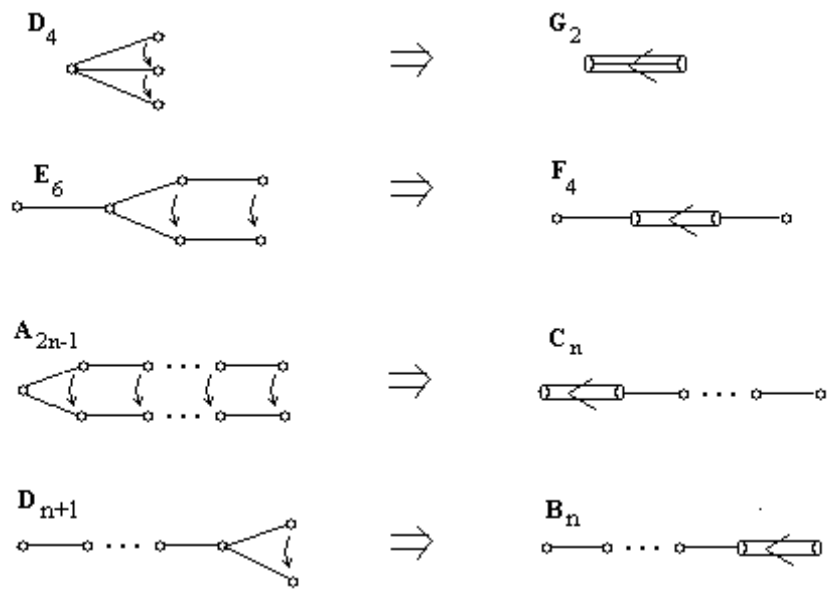}
\caption{\hspace{3mm} The folding operation applied to Dynkin
diagrams}
%%%%%% The label must come after caption
\label{folding}
\end{figure}

The \underline{folding of Dynkin diagrams} is defined by means of
the folding of the corresponding Cartan matrices. Let $\tau$ be a
diagram automorphism. The folded Cartan matrix $K^f$ is defined by
taking the sum over all $\tau$-orbits of the columns of $K$ (up to
some specific factor of this sum, Mohrdieck, \cite{Mohr04}).

\newpage
~\\ \vspace{2mm}
    \begin{color}{red}\LARGE\boldmath
      \section{\sc\bf The Slodowy correspondence, example: $\mathcal{T} \triangleleft \mathcal{O}$ }
      \end{color}
    \vspace{1mm}

\begin{figure}[h]
\centering
\includegraphics{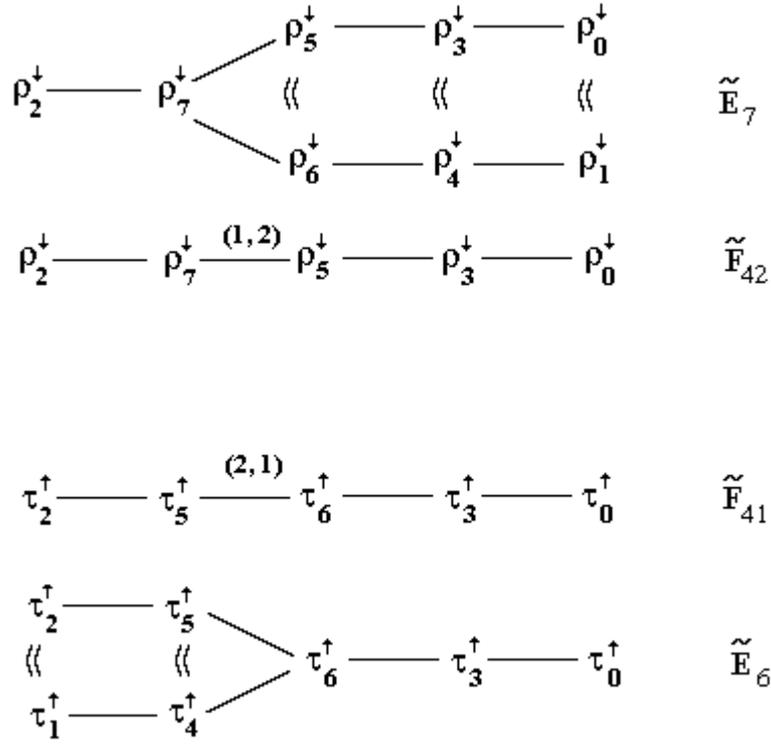}
\caption{\hspace{3mm} The induced and restricted representations
     of $\mathcal{T} \triangleleft \mathcal{O}$ }
%%%%%% The label must come after caption
\label{induced_restr_repr}
\end{figure}

We have

\begin{color}{blue}\Large\boldmath
\begin{equation*}
\label{tensor_prod_Slodowy_1}
\begin{split}
 & \tau_3\otimes\rho_0^\downarrow = \rho_3^\downarrow\otimes\rho_0^\downarrow =
\rho_3^\downarrow, \\
 & \tau_3\otimes\rho_2^\downarrow = \rho_3^\downarrow\otimes\rho_2^\downarrow =
\rho_7^\downarrow, \\
 & \tau_3\otimes\rho_3^\downarrow = \rho_3^\downarrow\otimes\rho_3^\downarrow =
\rho_0^\downarrow + \rho_5^\downarrow, \\
 & \tau_3\otimes\rho_5^\downarrow = \rho_3^\downarrow\otimes\rho_5^\downarrow =
\rho_3^\downarrow + \rho_7^\downarrow, \\
 & \tau_3\otimes\rho_7^\downarrow = \rho_3^\downarrow\otimes\rho_7^\downarrow =
\rho_2^\downarrow + 2\rho_5^\downarrow, \\
\end{split},
\hspace{5mm}
 \widetilde{A} = \left ( \begin{array}{ccccc}
        0 & 1 & 0 & 0 & 0 \vspace{2mm} \\
                1 & 0 & 2 & 0 & 0 \vspace{2mm} \\
                0 & 1 & 0 & 1 & 0 \vspace{2mm} \\
                0 & 0 & 1 & 0 & 1 \vspace{2mm} \\
                0 & 0 & 0 & 1 & 0
       \end{array}
  \right )
  \begin{array}{c}
    \rho_2^\downarrow \vspace{2mm} \\
    \rho_7^\downarrow \vspace{2mm} \\
    \rho_5^\downarrow \vspace{2mm} \\
    \rho_3^\downarrow \vspace{2mm} \\
    \rho_0^\downarrow
  \end{array}
\end{equation*}
\end{color}

\newpage
~\\ \vspace{2mm}
    \begin{color}{red}\LARGE\boldmath
      %%\section{\sc\bf The Kostant generating function }
      \section{\sc\bf Decomposition  $\pi_n |_G$ (Kostant) }
      \end{color}
    \vspace{1mm}

Let $\Sym(\mathbb{C}^2)$ be the symmetric algebra on $\mathbb{C}^2$,
in other words, $\Sym(\mathbb{C}^2) = \mathbb{C}[x_1, x_2]$. The
symmetric algebra $\Sym(\mathbb{C}^2)$ is a graded
$\mathbb{C}$-algebra:

\begin{color}{blue}\Large\boldmath
\begin{equation*}
    \Sym(\mathbb{C}^2) =
\mathop{\oplus}\limits_{m=0}^{\infty}{\Sym}^m(\mathbb{C}^2),
\end{equation*}
\end{color}\\
where ${\Sym}^m(\mathbb{C}^2)$ denotes the \underline{$m$th
symmetric power} of $\mathbb{C}^2$, which consists of the
homogeneous polynomials of degree $m$ in $x, y$:

\begin{color}{blue}\Large\boldmath
\begin{equation*}
    {\Sym}^m(\mathbb{C}^2) = \text{Span}\{x^m, x^{m-1}y, \dots, xy^{m-1}, y^m \}
\end{equation*}
\end{color}

For $n = 0,1,2,\dots$, let $\pi_n$ be the representation of $SU(2)$
in ${\Sym}^n(\mathbb{C}^2)$ induced by its action on $\mathbb{C}^2$.
The set \{$\pi_n \mid n \in \mathbb{Z}_+$\} is the set of all
irreducible representations of $SU(2)$.

Let $G$ be any finite subgroup of $SU(2)$. In \cite{Kos84}, Kostant
considered the following question:

\begin{color}{blue}\Large\boldmath
\medskip
{\sl How does $\pi_n |_G$ decompose for any $n \in \mathbb{N}$?}
\medskip
\end{color}

In other words: In the decomposition
\begin{color}{blue}\Large\boldmath
\begin{equation}
  \label{decomp_pi_n}
   \pi_n |_G = \sum\limits_{i=0}^r{m_i(n)\rho_i},
\end{equation}
\end{color}\\
where $\rho_i$ are irreducible representations of $G$, considered in
the context of the \begin{color}{blue}{\bf McKay
correspondence}\end{color},% Now, the question is the following one:

\begin{color}{blue}\Large\boldmath
\medskip
{\sl What are the multiplicities $m_i(n)$ equal to?}
\medskip
\end{color}

\newpage
    \begin{color}{red}\LARGE\boldmath
      %%\section{\sc\bf The Kostant generating function }
      \section{\sc\bf The Kostant generating function, the multiplicities $m_i(n)$  }
      \end{color}
    \vspace{1mm}

     In \cite{Kos84}, B.~Kostant obtained the multiplicities $m_i(n)$ by
studying the orbit structure of the Coxeter transformation on the
highest root of the corresponding root system.

The multiplicities $m_i(n)$ in (\ref{decomp_pi_n}) are calculated as
follows:

\begin{color}{blue}\Large\boldmath
\begin{equation*}
  \label{multipl_n}
     m_i(n) = \langle\pi_n | G , \rho_i\rangle.
\end{equation*}
\end{color}

 We extend the relation for multiplicity
 to the cases of \begin{color}{blue}{\bf restricted
 representations}\end{color}
  $\rho_i^\downarrow: = \rho_i\downarrow_H^G$
  and \begin{color}{blue}{\bf induced representations}\end{color}
   $\rho_i^\uparrow: = \rho_i\uparrow_H^G$, where $H$ is any subgroup of  $G$ %% \DL{?}
  (in the context of the Slodowy correspondence):

\begin{color}{blue}\Large\boldmath
\begin{equation*}
  \label{decomp_pi_n_2}
   m_i^\downarrow(n) = \langle\pi_n | H ,\rho_i^\downarrow\rangle,   \hspace{3mm}
   m_i^\uparrow(n) = \langle\pi_n | G , \rho_i^\uparrow\rangle.
\end{equation*}
\end{color}

Kostant introduced the \begin{color}{blue}{\bf generating function}
\boldmath$P_G(t)$\end{color} as follows:

\begin{color}{blue}\Large\boldmath
\begin{equation}
 \label{Kostant_gen_func}
   P_G(t)  =
     \left(
       \begin{array}{cc}
         & [P_G(t)]_0 \vspace{2mm} \\
         & \ldots \vspace{2mm} \\
         & [P_G(t)]_r \vspace{2mm}
       \end{array} \right )
           :=
     \left(
       \begin{array}{cc}
         & \sum\limits_{n=0}^{\infty}m_0(n){t^n} \\
         & \ldots  \\
         & \sum\limits_{n=0}^{\infty}m_r(n){t^n}
       \end{array} \right ).
\end{equation}
\end{color}

We introduce
\begin{color}{blue}\boldmath$P_{G\uparrow}(t)$\end{color} (resp.
\begin{color}{blue}\boldmath$P_{G\downarrow}(t)$\end{color}) by substituting  $m_i^\uparrow(n)$ (resp.
$m_i^\downarrow(n)$) instead of $m_i(n)$.

\begin{color}{blue}\large\boldmath
\begin{equation}
 \label{Kostant_gen_func_2}
   P_{G\uparrow}(t)   :=
     \left(
       \begin{array}{cc}
         & \sum\limits_{n=0}^{\infty}m_0^\uparrow(n){t^n} \\
         & \ldots  \\
         & \sum\limits_{n=0}^{\infty}m_r^\uparrow(n){t^n}
       \end{array} \right ),
       \quad
   P_{G\downarrow}(t)   :=
     \left(
       \begin{array}{cc}
         & \sum\limits_{n=0}^{\infty}m_0^\downarrow(n){t^n} \\
         & \ldots  \\
         & \sum\limits_{n=0}^{\infty}m_r^\downarrow(n){t^n}
       \end{array} \right ).
\end{equation}
\end{color}

\newpage
~\\ \vspace{2mm}
    \begin{color}{red}\LARGE\boldmath
      \section{\sc\bf The Poincar\'{e} series for the binary polyhedral groups}
      \end{color}
    \vspace{1mm}

The multiplicity $m_0(n)$ corresponds to the trivial representation
$\rho_0$ in ${\Sym}^n(\mathbb{C}^2)$. The
\begin{color}{blue}{\bf algebra of invariants}
\boldmath$R^G$\end{color} coincides with $\Sym(\mathbb{C}^2)$, and
$[P_G(t)]_0$ is the Poincar\'{e} series of the algebra of invariants
$R^G = \Sym(\mathbb{C}^2)^G$, i.e., %% \DL{protivorechie: prochti frazu}
(Kostant, \cite{Kos84})
%% \DL{ne opredelen ryad $P(A,t)$, gde $A$ --
%% kakayato (veroyatno, Z-graduirovannaya s eshche paroj svojstv, algebra)}

\begin{color}{blue}\boldmath
\begin{equation*}
  \label{poincare_alg_inv}
    [P_G(t)]_0 = P(\Sym(\mathbb{C}^2)^G,t).
\end{equation*}
\end{color}

  .\\

\begin{theorem} [Kostant, Kn\"{o}rrer, Gonzalez-Sprinberg, Verdier]
  The Poincar\'{e} series $[P_G(t)]_0$ can be calculated as the
following rational function:

\begin{color}{blue}\boldmath
 \begin{equation*}
  \label{K_K_GV}
      [P_G(t)]_0 = \frac{1 + t^h}{(1 - t^a)(1 - t^b)},
\end{equation*}
 \end{color} \\
where $h$ is the Coxeter number, while $a$ and $b$ are given by the
system
%% \DL{chto vsegda est' celye resheniya? eto, chto-li,
%% ochevidno?! mozhet tak i skazat'...}

\begin{color}{blue}\boldmath
 \begin{equation*}
  \label{Kostant_numbers_a_b}
    a+  b = h + 2, \quad ab = 2|G|.
 \end{equation*}
\end{color}

\end{theorem}

\newpage
    \begin{color}{red}\LARGE\boldmath
      \section{\sc\bf The McKay-Slodowy operator }
      \end{color}
    \vspace{1mm}

 We set
\begin{color}{blue}\large\boldmath
 \begin{equation*}
   %%\label{Kostant_relation}
       v_n =
      \begin{cases}
        & \sum\limits_{i=0}^r{m_i(n)}\alpha_i, \quad  \text { for } B = A,  \\
        & \sum\limits_{i=0}^r{m_i^\downarrow(n)}\alpha_i, \quad \text { for } B = \widetilde{A}, \\
        & \sum\limits_{i=0}^r{m_i^\uparrow(n)}\alpha_i, \quad \text { for } B = \widetilde{A}^\vee, \\
      \end{cases},
  \quad
   \rho_i =
     \begin{cases}
        \rho_i \\
        \rho_i^\downarrow  \\
        \rho_i^\uparrow  \\
     \end{cases},
 \quad
    m_i(n) =
     \begin{cases}
        m_i(n)  \\
        m_i^\downarrow(n)  \\
        m_i^\uparrow(n) \\
     \end{cases}.
 \end{equation*}
\end{color}

The following result of B.~Kostant \cite{Kos84}, which holds for the
McKay operator (\ref{main_McKay}) holds also for the Slodowy
operators (\ref{matrix_Slodowy}), (\ref{matrix_Slodowy_1}).
\vspace{2mm} \\

\begin{proposition}
  \label{kostant_prop}
  If $B$ is either the McKay operator $A$ or one of the Slodowy
  operators
  $\widetilde{A}$ or $\widetilde{A}^\vee$, then

\begin{color}{blue}\Large \boldmath
 \begin{equation}
   \label{Kostant_relation}
     Bv_n = v_{n-1} + v_{n+1}.
 \end{equation}
\end{color}
\end{proposition}

\PerfProof  We have
\begin{color}{blue}\large\boldmath
\begin{equation*}
    Bv_n = B
    \left (
    \begin{array}{c}
       m_0(n) \\
       \dots  \\
       m_r(n)
    \end{array}
    \right ) =     \left (
    \begin{array}{c}
       \sum a_{0i}\langle\rho_i, \pi_n\rangle \\
       \dots  \\
       \sum a_{ri}\langle\rho_i, \pi_n\rangle \\
    \end{array}
    \right ) =
    \left (
    \begin{array}{c}
       \langle\rho\otimes\rho_0, \pi_n\rangle \\
       \dots  \\
       \langle\rho\otimes\rho_r, \pi_n\rangle \\
    \end{array}
    \right ),
\end{equation*}
\end{color} \\
where $\rho$ is the irreducible $2D$ representation which coincides
with the representation $\pi_1$ in ${\Sym}^2(\mathbb{C}^2)$. For
representations $\rho_i$ of any finite subgroup $G \subset SU(2)$,
we have
\begin{color}{blue}\large\boldmath$\langle\chi_i\chi_j, \chi_k\rangle ~=~
\langle\chi_i, \chi_j\chi_k\rangle$\end{color}, and

\begin{color}{blue}\large\boldmath
\begin{equation*}
  \label{McKay_oper_2}
   Bv_n =
    \left (
    \begin{array}{c}
       \langle\pi_1\otimes\rho_0, \pi_n\rangle \\
       \dots  \\
       \langle\pi_1\otimes\rho_r, \pi_n\rangle \\
    \end{array}
    \right ) =
     \left (
    \begin{array}{c}
       \langle\rho_0, \pi_1\otimes\pi_n\rangle \\
       \dots  \\
       \langle\rho_r, \pi_1\otimes\pi_n\rangle \\
    \end{array}
    \right ).
\end{equation*}
\end{color} \\

By Clebsch-Gordan formula we have

\begin{color}{blue}\large\boldmath
\begin{equation*}
  \label{McKay_oper_4}
   \pi_1\otimes\pi_n = \pi_{n-1} \oplus \pi_{n+1},
\end{equation*}
\end{color} \\
where $\pi_{-1}$ is the zero
%% \DL{eto ved' ne "trivial". Opredelit'
%%by, raz uzh eto populyarnyj text i dazhe trivial rep opredelen ...}
representation. \qed

\newpage
~\\ \vspace{2mm}
    \begin{color}{red}\LARGE\boldmath
      \section{\sc\bf The McKay-Slodowy operator (2) }
      \end{color}
    \vspace{1mm}

Let $x = \widetilde{P}_G(t)$ be given by (\ref{Kostant_gen_func}),
(\ref{Kostant_gen_func_2}), namely:

\begin{color}{blue}\Large\boldmath
\begin{equation}
  \label{gen_poincare_ser_3}
    \widetilde{P}_G(t) =
     \begin{cases}
        & P_G(t) \hspace{3mm}\text{ for } B = A, \\
        & \\
        & P_{G\downarrow}(t) \hspace{2mm}\text{ for }  B = \widetilde{A}, \\
        & \\
        & P_{G\uparrow}(t) \hspace{2mm}\text{ for }  B = \widetilde{A}^\vee,\\
    \end{cases}
\end{equation}
\end{color}

\begin{proposition} We have
\begin{color}{blue}\Large\boldmath
\begin{equation}
  \label{McKay_oper_5}
  tBx = (1 + t^2)x - v_0,
\end{equation}
\end{color} \\
where $B$ is either the McKay operator $A$ or one of the Slodowy
operators
  $\widetilde{A}$, $\widetilde{A}^\vee$.
\end{proposition}

\PerfProof
  From (\ref{Kostant_relation}) we obtain

\begin{color}{blue}\Large\boldmath
\begin{equation*}
\begin{split}
 Bx = & \sum\limits_{n=0}^{\infty}Bv_n{t^n} =
        \sum\limits_{n=0}^{\infty}(v_{n-1} + v_{n+1}){t^n} = \\
 &       \sum\limits_{n=0}^{\infty}v_{n-1}t^n +
        \sum\limits_{n=0}^{\infty}v_{n+1}t^n = \\
 &   t\sum\limits_{n=1}^{\infty}v_{n-1}t^{n-1} +
    t^{-1}\sum\limits_{n=0}^{\infty}v_{n+1}t^{n+1} = \\
 &  tx + t^{-1}(\sum\limits_{n=0}^{\infty}v_n{t^n} - v_0) =
    tx + t^{-1}x - t^{-1}v_0.  \qed
\end{split}
\end{equation*}
\end{color} \\

\newpage
~\\ \vspace{2mm}
    \begin{color}{red}\LARGE\boldmath
      \section{\sc\bf The Ebeling theorem }
        \label{Ebeling_th}
      \end{color}
    \vspace{1mm}

 W.~Ebeling in \cite{Ebl02} established the connection between the Poincar\'{e}
 series, the Coxeter transformation {\bf C}, and the corresponding affine
 Coxeter transformation  ${\bf C}_a$ (in the context of the McKay correspondence).
\vspace{2mm}

\begin{theorem}
 {\it  Let $G$ be a binary polyhedral group and let $[{P}_G(t)]_0$
  be the Poincar\'{e} series. Then

\begin{color}{blue}\Large\boldmath
\begin{equation*}
       [{P}_G(t)]_0 = \frac{\det{M_0}(t)}{\det{M}(t)},
\end{equation*}
where
\begin{equation*}
   \det{M}(t) = \det|t^{2}I - {\bf C}_a|, \hspace{5mm}
   \det{M_0}(t) = \det|t^{2}I - {\bf C}|,
\end{equation*}
\end{color} \vspace{2mm}\\
{\bf C} is the Coxeter transformation and ${\bf C}_a$ is the
corresponding affine Coxeter transformation}. \vspace{2mm}
\end{theorem}

    We extend this fact to the case of multiply-laced diagrams,
    and generalized Poincar\'{e} series $[\widetilde{P}_G(t)]_0$
    (in the context of the McKay-Slodowy correspondence), namely:

\begin{color}{blue}\Large\boldmath
\begin{equation}
  \label{gen_Webiling_th}
       [\widetilde{P}_G(t)]_0 = \frac{\det{M_0}(t)}{\det{M}(t)},
\end{equation}
\end{color} \vspace{2mm}\\
 see (\ref{gen_poincare_ser_3}).

 \newpage
    \begin{color}{red}\LARGE\boldmath
      \section{\sc\bf The Ebeling theorem (2) }
      \end{color}
    \vspace{2mm}

 {\it Proof} of (\ref{gen_Webiling_th}).  From (\ref{McKay_oper_5})
 we have

\begin{color}{blue}\Large\boldmath
  \begin{equation*}
   [(1 + t^2)I - tB]x = v_0,
\end{equation*}
\end{color} \\
where $x$ is the vector $\widetilde{P}_G(t)$ and by Cramer's rule
the first coordinate of $\widetilde{P}_G(t)$ is

\begin{color}{blue}\Large\boldmath
 \begin{equation*}
       [\widetilde{P}_G(t)]_0 = \frac{\det{M_0}(t)}{\det{M}(t)},
 \end{equation*}
\end{color} \\
 where

 \begin{color}{blue}\Large\boldmath
 \begin{equation*}
  \label{M_def_affine}
       \det{M}(t) = \det\left((1 + t^2)I - tB\right),
 \end{equation*}
 \end{color} \\
and $M_0(t)$ is the matrix obtained by replacing the first column of
$M(t)$ by $v_0 = (1,0,...,0)^t$. The vector $v_0$ corresponds to the
trivial representation $\pi_0$, and by the McKay-Slodowy
correspondence, $v_0$ corresponds to the particular vertex which
extends the Dynkin diagram to the extended Dynkin diagram.
 (For calculation of $v_0$, see (\ref{frob_2}), (\ref{clifford_2}), and
 Remark \ref{rem_v0}). Therefore, if $\det{M}(t)$ corresponds to the affine Coxeter
transformation, and

 \begin{color}{blue}\Large\boldmath
 \begin{equation}
 \label{M_affine_cox}
       \det{M}(t) = \det|t^{2}I - {\bf C}_a|,
 \end{equation}
 \end{color} \\
then $\det{M_0}(t)$ corresponds to the Coxeter transformation, and

 \begin{color}{blue}\Large\boldmath
 \begin{equation*}
 \label{M_cox}
       \det{M_0}(t) = \det|t^{2}I - {\bf C}|.
 \end{equation*}
 \end{color} \\
So, it suffices to prove (\ref{M_affine_cox}), i.e.,

 \begin{color}{blue}\Large\boldmath
 \begin{equation}
 \label{M_affine_cox_2}
        \det[(1 + t^2)I - tB] = \det|t^{2}I - {\bf C}_a|.
 \end{equation}
 \end{color}

 \newpage
    \begin{color}{red}\LARGE\boldmath
      \section{\sc\bf The Ebeling theorem (3) }
      \end{color}
    \vspace{2mm}

If $B$ is the McKay operator $A$ given by (\ref{main_McKay}), then

 \begin{color}{blue}\Large\boldmath
\begin{equation*}
 \label{McKay_operator}
       B = 2I - K =
    \left ( \begin{array}{cc}
            0     & -2D \\
            -2D^t &  0
            \end{array}
    \right ),
 \end{equation*}
 \end{color} \\
where $K$ is a symmetric Cartan matrix (\ref{symmetric_B}). If $B$
is the Slodowy operator $\widetilde{A}$ or $\widetilde{A}^\vee$
given by (\ref{lab_slodowy_1}), (\ref{lab_slodowy_2}), then

%% \DL{ya pravda ne ponimayu: a kakie mogut byt' vyrozhdeniya? po
%% kakomu parametru?}
%% R.S. Ne ponyal voprosa: to v.2
 \begin{color}{blue}\Large\boldmath
 \begin{equation*}
 \label{Slodowy_operator}
       B = 2I - K =
    \left ( \begin{array}{cc}
            0     & -2D \\
            -2F &  0
            \end{array}
    \right ),
 \end{equation*}
 \end{color} \\
where $K$ is the symmetrizable Cartan matrix (\ref{matrix_K}). Thus,
in the generic case

 \begin{color}{blue}\Large\boldmath
 \begin{equation}
 \label{M_t}
       M(t) = (1+t^2)I - tB =
    \left ( \begin{array}{cc}
            1+t^2 &  2tD     \\
            2tF  &  1 + t^2
            \end{array}
     \right ).
 \end{equation}
 \end{color} \\

Assuming $t \neq 0$ we deduce from (\ref{M_t}) that
 \begin{color}{blue}\Large\boldmath
\begin{equation}
 \label{M_t_2}
 \begin{array}{cc}
       M(t)
    \left ( \begin{array}{c}
            x    \\
            y
            \end{array}
    \right ) = 0 & \Longleftrightarrow
    \left \{
     \begin{array}{c}
            (1 + t^2)x = -2tDy,  \\
            2tFx = -(1 + t^2)y.
            \end{array}
     \right .  \vspace{5mm} \\
    & \Longleftrightarrow
    \left \{
     \begin{array}{c}
            \displaystyle\frac{(1 + t^2)^2}{4t^2}x = FDy,   \vspace{3mm} \\
            \displaystyle\frac{(1 + t^2)^2}{4t^2}y = DFy.
            \end{array}
     \right .
  \end{array}
\end{equation}
\end{color}

 According to (\ref{DDt_DtD}), and the propositions about Jordan  normal form of the Coxeter
 transformation,  we see that $t^2$ is an eigenvalue
 of the affine Coxeter transformation ${\bf C}_a$, i.e., (\ref{M_affine_cox_2}) together
with (\ref{M_affine_cox}) are proved. \qedsymbol

  For further details and references,  see \cite{St08}.
  For applications to the singularity theory, see \cite{Ebl08}.

 \newpage
    \begin{color}{red}\LARGE\boldmath
      \section{\sc\bf Proportionality of characteristic polynomials and folding}
      \end{color}
    \vspace{2mm}

   By calculating, we obtain that Poincar\'{e} series coincide for the following pairs
   of diagrams
\begin{color}{blue}\Large\boldmath
 \begin{equation*}
   \begin{array}{cc}
       {D}_4  \text{ and } {G}_2, \hspace{7mm}
      & {E}_6  \text{ and } {F}_4, \\
       {D}_{n+1} \text{ and } {B}_n  (n \geq 4), \hspace{7mm}
      &  {A}_{2n-1} \text{ and } {C}_n.
   \end{array}
\end{equation*}
\end{color} \\
Note that the second elements of the pairs are obtained by
\underline{folding}:

\begin{color}{blue}\Large\boldmath
 \begin{equation*}
   \frac{\mathcal{X}({D}_4)}{{\mathcal{X}(\widetilde{D}_4)}} =
   \frac{\mathcal{X}({G}_2)}{{\mathcal{X}(\widetilde{G}_{21})}} =
    \frac{\lambda^3 + 1}{(\lambda^2 - 1)^2}.
\end{equation*}
\end{color} \\

\begin{color}{blue}\Large\boldmath
 \begin{equation*}
   \frac{\mathcal{X}({E}_6)}{{\mathcal{X}(\widetilde{E}_6)}} =
   \frac{\mathcal{X}({F}_4)}{{\mathcal{X}(\widetilde{F}_{41})}} =
    \frac{\lambda^6 + 1}{(\lambda^4 -  1)(\lambda^3 - 1)}.
\end{equation*}
\end{color} \\

\begin{color}{blue}\Large\boldmath
 \begin{equation*}
   \frac{\mathcal{X}({D}_{n+1})}{{\mathcal{X}(\widetilde{D}_{n+1})}} =
   \frac{\mathcal{X}({B}_n)}{{\mathcal{X}(\widetilde{B}_{n})}} =
    \frac{\lambda^{n} + 1}{(\lambda^{n-1} -  1)(\lambda^2 - 1)}.
\end{equation*}
\end{color} \\

\begin{color}{blue}\Large\boldmath
 \begin{equation*}
   \frac{\mathcal{X}({A}_{2n-1})}{{\mathcal{X}(\widetilde{A}_{2n-1})}} =
   \frac{\mathcal{X}({C}_n)}{{\mathcal{X}(\widetilde{C}_{n})}} =
    \frac{\lambda^{n} + 1}{(\lambda^n -  1)(\lambda - 1)}.
\end{equation*}
\end{color} \\

\newpage
    \begin{color}{red}\LARGE\boldmath
      \section{\sc\bf Acknowledgements }
      \end{color}
    \vspace{12mm}

\begin{color}{blue}\LARGE
I would like to thank the organizers of the Workshop {\bf Jose
Antonio de la Pena, Vlastimil Dlab,} and {\bf Helmut Lenzing} who
gave me the opportunity to present this survey.
\end{color}

\begin{color}{blue}\LARGE
I am thankful to {\bf John McKay} and {\bf Dimitry Leites} for
 helpful comments  to this survey.
\end{color}

 \vspace{3mm}

\newpage
\begin{color}{blue}

~\\ 
\end{color}

\vspace{3mm}

\end{document}